\documentclass[final,12pt]{colt2026} 

\title[Score-based modular reduction]{
Fast Score-Based Sampling via Log-Concave Reductions}

\usepackage{times}

\coltauthor{%
  \Name{Martin J. Wainwright}  \\
  \Email{mjwain@mit.edu} \\
 \addr  Laboratory for Information and Decision Systems  \\
 Statistics and Data Science Center \\
 EECS \& Mathematics \\
 Massachusetts Institute of Technology
}

\usepackage{url}
\usepackage{mathtools}
\usepackage{mathbbol}
\usepackage{mathrsfs}  
\usepackage{multicol}
\usepackage{rotating}
\usepackage{enumitem, comment, xifthen}

\usepackage{xcolor}
\usepackage{booktabs}
\usepackage[mathscr]{euscript}
\usepackage[most]{tcolorbox}

\usepackage{natbib}

\usepackage{aliascnt}

\newtheorem{lemma}{Lemma}
\newtheorem{proposition}{Proposition}
\newtheorem{corollary}{Corollary}

\usepackage[nameinlink]{cleveref}

\crefname{theorem}{Theorem}{Theorems}
\Crefname{theorem}{Theorem}{Theorems}

\crefname{lemma}{Lemma}{Lemmas}
\Crefname{lemma}{Lemma}{Lemmas}

\crefname{proposition}{Proposition}{Propositions}
\Crefname{proposition}{Proposition}{Propositions}

\crefname{corollary}{Corollary}{Corollaries}
\Crefname{corollary}{Corollary}{Corollaries}

\hypersetup{
  colorlinks,
  linkcolor={red!72!black},
  citecolor={blue!72!black},
  urlcolor={magenta!72!black}
}

\newcommand\blfootnote[1]{%
  \begingroup
  \renewcommand\thefootnote{}\footnote{#1}%
  \addtocounter{footnote}{-1}%
  \endgroup
}

\usepackage{rat_macros}
\usepackage{tearless_macros}

\newcommand{\Mtar}{\ensuremath{L}}
\newcommand{\Cmax}{\ensuremath{\CovCon_{\max}}}
\newcommand{\utarget}{u^\star}
\newcommand{\Kfinal}{\ensuremath{K}}
\newcommand{\vup}[1]{\ensuremath{v^{(#1)}}}
\newcommand{\delup}[1]{\ensuremath{\Delta^{(#1)}}}
\newcommand{\MyInt}[1]{\ensuremath{\mathcal{I}^{(#1)}}}
\newcommand{\Nint}[1]{\ensuremath{N^{(#1)}}}
\newcommand{\var}{\ensuremath{\operatorname{var}}}

\newcommand{\Markov}{\ensuremath{\mathcal{P}}}
\newcommand{\Wass}{\ensuremath{\mathscr{W}}}
\newcommand{\qtil}{\ensuremath{\widetilde{q}}}
\newcommand{\vtil}{\ensuremath{\widetilde{v}}}
\newcommand{\thetatil}{\ensuremath{\widetilde{\theta}}}
\newcommand{\thetabar}{\ensuremath{\bar{\theta}}}
\newcommand{\Kull}{D_{\mbox{\tiny{KL}}}}

\newcommand{\Nit}{N}
\newcommand{\NitSLC}{\ensuremath{\subtiny{\Nit}{\texttt{SLC}}}}

\newcommand{\Nula}{\ensuremath{\subtiny{\Nit}{\texttt{ULA}}}}

\newcommand{\Nfast}{\ensuremath{\subtiny{\Nit}{\texttt{fast}}}}

\newcommand{\dtv}{\ensuremath{\subsmall{D}{\mbox{TV}}}}
\newcommand{\dtvsquared}{D^2_{\scaleto{\mbox{TV}}{5pt}}}
\newcommand{\dkl}{\ensuremath{\subsmall{D}{\mbox{KL}}}}

\newcommand{\pvcondu}{p_{\myspecial{v} \mid \myspecial{u}}}
\newcommand{\pucondv}{p_{\myspecial{u} \mid \myspecial{v}}}

\newcommand{\pkback}{\ensuremath{p_{k \mid k+1}}}

\newcommand{\mlowerbar}{\mlower}

\newcommand{\CovCon}{\ensuremath{B}}
\newcommand{\IterB}{\ensuremath{\lambda}}

\newcommand{\ptarget}{\ensuremath{{p_{\myspecial{x}}}}}

\usepackage{xspace}
\usepackage{enumitem}

\newcommand{\Otil}{\ensuremath{\widetilde{O}}}

\newcommand{\kull}{\ensuremath{D_{\scaleto{\text{KL}}{5pt}}}}

\newcommand{\Tfinal}{\ensuremath{T}}

\newcommand{\hackalpha}{\alpha}
\newcommand{\hackbeta}{\beta}
\newcommand{\Tworst}{T_{\scaleto{worst}{5pt}}}
\newcommand{\polylog}{\ensuremath{\operatorname{polylog}}}

\newcommand{\MyProb}{\ensuremath{\mathscr{S}}}

\newcommand{\Back}[1]{\ensuremath{#1}}
\newcommand{\Qback}{\Back{\mathcal{Q}}}
\newcommand{\Qbacktil}{\Back{\widetilde{\mathcal{Q}}}}
\newcommand{\Pback}{\Back{\mathcal{P}}}
\newcommand{\scorehat}[1]{\widehat{s}_{#1}}

\newcommand{\qbar}{\ensuremath{\bar{q}}}
\renewcommand{\scon}{\alpha}

\newcommand{\pihat}{\widehat{\pi}}

\begin{document}

\maketitle

\blfootnote{Accepted for presentation at the Conference on Learning Theory (COLT) 2026.}

\begin{abstract}%
Sampling based on score diffusions has led to striking empirical
results, and has attracted considerable attention from various
research communities.  It depends on availability of (approximate)
Stein score functions for various levels of additive noise.  We show
how in some generality, the availability of scores allows the general
problem to be ``reduced'' to sampling from an adaptively constructed
sequence of $K$ strongly log-concave (SLC) sub-problems.  The
reduction is simple, constructive and algorithm-independent, so that
any SLC sampler can be used as a subroutine.  Various bounds on
score-based sampling complexity follow directly: for instance,
high-accuracy SLC samplers yield $\Otil(K \sqrt{d}
\polylog(1/\varepsilon))$ guarantees for accuracy $\varepsilon$ in
dimension $d$, where randomized midpoint SLC schemes yield $\Otil(K 
d^{1/3} \operatorname{poly}(1/\varepsilon))$ guarantees.  When the
original distribution itself is SLC, we prove that $K \leq 1 +
\log_2(\kappa)$, thereby obtaining the first efficient procedure with
logarithmic dependence on condition number $\kappa$; for general
distributions, the quantity $K$ depends on the geometry of score
Hessian across the trajectory.  Our analysis is direct and simple,
involving techniques and insights complementary to those in standard
analyses of discretized diffusions.
\end{abstract}


\section{Introduction}

The problem of drawing samples from a $d$-dimensional density is a
core computational challenge.  Efficient samplers are essential for
Monte Carlo approximation
(e.g.,~\cite{robert2004monte,rubinstein2008simulation}); exploration
of posterior distributions in Bayesian statistics and inverse problems
(e.g.,~\cite{gelman2013bayesian,brooks2011handbook}); and generation
of images, audio and other structured data in generative AI
(e.g.,~\cite{rombach2022high,croitoru2023diffusion_vision,chen2024overview}).

\paragraph{Score-based diffusions:}
In recent years, researchers have demonstrated dramatic advances in
sampling through the use of score-based diffusion
models~\citep{sohldickstein2015deep,song2019generative,ho2020denoising,song2021score}.
All these procedures are based on a forward noising process: beginning
with a sample $X$ from the target distribution $\ptarget$, it converts
it to some form of ``noise'', most often a standard Gaussian vector.
In continuous time, this forward process can be described by a
stochastic differential equation (SDE), and the problem of drawing
samples corresponds to simulating the evolution of the reverse-time
SDE~\citep{haussmann1986time,anderson1982reverse} that tracks backward
from the noise $W$ to a fresh sample from $\ptarget$.  Stochastic
sampling schemes are based on careful discretizations of this
reverse-time SDE (e.g.,~\cite{ho2020denoising,song2021score,
  lee2022convergence,li2023towards,lee2023convergence,chen2023sampling,chen2023improved,benton2024nearly}),
whereas other sampling schemes make use of an ordinary differential
equation (ODE) that describes the backwards evolution
(e.g.,~\cite{song2021score,benton2023error,albergo2023stochastic,chen2023probability,cai2025minimax}).
In both cases, the forward process is useful, it provides the data
needed to estimate the Stein score functions that describe the
backwards evolution, using methods such as score matching or
Tweedie-based denoising
(e.g.,~\cite{robbins1956empirical,Miyasawa1961,hyvarinen2005estimation,vincent2011connection}).
There are now a wide variety of schemes within the general diffusion
framework along with a relatively rich theoretical understanding;
see~\Cref{SecRelated} for further discussion.

\paragraph{Fast sampling for ``nice'' distributions:}
In parallel, the past ten years have witnessed tremendous advances in
the problem of drawing samples from ``nice'' distributions, based only
on first-order information from the original distribution (and not
invoking a diffusion path).  For instance, there are highly efficient
methods, along with associated theoretical guarantees, for sampling
from strongly log-concave (SLC) distributions, as well as those that
satisfy a geometric inequality (e.g., log-Sobolev, or Poincar\'{e}).
A wide spectrum of methods have been studied, including the unadjusted
Langevin algorithm (ULA), its Metropolis-corrected variant (MALA),
higher-order extensions including Hamiltonian Monte Carlo, as well as
various proximal schemes
(e.g.,~\cite{Dal16,cheng2018convergence,durmus2017nonasymptotic,DwiCheWaiYu19,MouFlaWaiBar19b,Che+19_HMC,vempala_wibisono_2022,chewi2022analysis,chen_gatmiry_2023}).
The modular scheme of this paper allows \emph{any strongly log-concave
(SLC) sampler} to be applied, and we obtain a spectrum of rates
depending on this choice.  Notably, for SLC problems with order-one
conditioning, there exist high-accuracy Metropolized samplers with
$\sqrt{d} \: \polylog(1/\varepsilon)$ iteration
complexity~\citep{altschuler_chewi_2023,Chewi_book}.  While the
standard Metropolis implementation requires density evaluations, we
note that very recent work~\citep{Sasha} provides a randomized scheme
that achieves the Metropolized complexity using only first-order
gradient information.

\paragraph{Our contributions:} In this paper, we bring these
two lines of research into close contact, in particular by describing
and analyzing a simple and modular scheme for score-based sampling.
We show how, given the availability of annealed Stein scores, it is
possible to ``reduce'' the problem of sampling from a general target
density $\ptarget$ to a sequence of $K$ calls to \emph{any SLC
sampler} applied to distributions with condition number at most $2$.
We give guarantees in which number of calls $K$ to the SLC sampler is
independent of the target accuracy $\varepsilon$, but depends on the
geometric structure of the problem.

\noindent More specifically, we prove two main theorems:
\vspace*{-0.09in}
\begin{itemize}[leftmargin=*, itemsep=0pt]
\item In~\Cref{ThmConcave}, we study the problem of sampling from a
  strongly log-concave density on $\real^d$ with condition number
  $\kappa$.  This is a very well-studied problem, and standard
  Langevin-type procedures exhibit \emph{linear} scaling with
  $\kappa$.  By adapting our modular scheme that exploits Stein
  scores, we show that this dependence can be reduced to
  \emph{logarithmic} in $\kappa$.  In particular, we do so by reducing
  to $K+1$ SLC problems with $K \leq 1 + \log_2(\kappa)$. To the best
  of our knowledge, this is the first efficient scheme that exhibits
  such logarithmic dependence.  When high-accuracy SLC samplers are
  used within the reduction, we guarantee an overall iteration
  complexity of $\Otil \big( \sqrt{d} \; \log(\kappa)
  \polylog(1/\varepsilon) \big)$.

\item In~\Cref{ThmMultiModal}, we study the use of our modular scheme
  for sampling from a general multi-modal density.  We specify a
  forward trajectory of length $K$, specified by an \emph{adaptive
  stepsize sequence}, that ensures a SLC sequence of sub-problems.
  This leads to concrete bounds on sampling complexity in terms of
  trajectory length $K$, and the complexity of solving these
  sub-problems, and we exhibit a sampler with iteration complexity
  $\Otil( K \sqrt{d} \polylog(1/\varepsilon))$.  In
  Corollary~\ref{CorMultiWorst}, we provide a worst-case bound on $K$
  in terms of a geometric Lipschitz constant, but suspect that this
  guarantee can be improved.
\end{itemize}
We put these results in the context of past work in the next section,
as well as in the discussion following the statement of our theorems.


\subsection{Related work}
\label{SecRelated}

There is a long line of work on fast algorithms for sampling from
strongly log-concave distributions (SLC), as well as more general
families, including those satisfying log-Sobolev and Poincar\'{e}
inequalities (e.g.,~\cite{Dal16,
  cheng2018convergence,durmus2017nonasymptotic,DwiCheWaiYu19,MouFlaWaiBar19b,Che+19_HMC,vempala_wibisono_2022,chewi2022analysis,chen_gatmiry_2023}).
The modular reduction in this paper allows any of these procedures to
be called as a black box routine.  Various algorithms have been
analyzed, including the unadjusted Langevin (ULA) algorithm, its
Metropolis-adjusted variant (known as MALA), higher-order schemes
including randomized midpoint and Hamiltonian Monte Carlo; and
samplers based on proximal updates.  Suitable variants can achieve
iteration complexity proportional to $\sqrt{d}$; of particular
relevance to this paper are high-accuracy samplers for SLC
distributions with iteration complexity scaling polynomially in
$\log(1/\varepsilon)$
(e.g.,~\cite{DwiCheWaiYu19,Che+19_HMC,chen_gatmiry_2023,altschuler_chewi_2023}).
Other schemes, such as randomized midpoint
discretizations~\citep{ShenLee2019RandomizedMidpoint}, sacrifice the
$\polylog(1/\varepsilon)$ scaling but reduce dimension dependence to
$d^{1/3}$.  Our general scheme allows any of these SLC samplers to be
applied.

For diffusion-based samplers, there is now a wide range of theoretical
results, applying to both stochastic (SDE-based) samplers
(e.g.,~\cite{lee2022convergence,li2023towards,lee2023convergence,chen2023sampling,chen2023improved,benton2024nearly})
as well as (ODE or flow-based) deterministic ones
(e.g.,~\cite{song2021score,albergo2023stochastic,benton2023error,chen2023probability,cai2025minimax}).
Earlier analyses of the iteration complexity, meaning the number of
iterations needed to obtain $\varepsilon$-accurate samples, exhibited
polynomial scaling in the dimension.  Focusing on the KL divergence,
recent results have reduced this dependence to linear in dimension $d$
for both stochastic
samplers~\citep{conforti2023score,benton2024nearly} and ODE-based
samplers~\citep{li_wei_chi_chen_2024}, and both classes of methods
have polynomial scaling in $(1/\varepsilon)$.  By comparison, our
modular scheme yields a method with KL iteration complexity scaling as
$\sqrt{d}$, and polynomially in $\log(1/\varepsilon)$;
see~\Cref{ThmMultiModal} for details.

The idea of decomposing a ``hard'' sampling problem into strongly
log-concave problems can be understood as an auxiliary variable
method~\citep{Liu01,Higdon98}, and has been exploited for sampling
from diffusions~\citep{ZhangNew}, as well as posterior distributions
in neural networks~\citep{curtis2024}, and sparse linear
regression~\citep{MonWu24}.  Mostly closely related to our work is the
paper of~\citet{ZhangNew}, who proposed two different algorithms
(RTK-ULA and RTK-MALA) that exploit log-concave segments within a
diffusion.  Among other results, they showed that the RTK-MALA
procedure has TV mixing time scaling as $d^2 \polylog(1/\varepsilon)$,
which is the first high-accuracy guarantee for a diffusion-based
procedure.  Our work extends this idea of extracting SLC sub-problems
in a way that is \emph{adaptive}, depending on the problem structure,
and \emph{algorithm-independent}, thereby allowing any SLC sampler to
be applied to the sub-problems.  Consequently, our guarantees inherit
the best available bounds for the SLC subroutine for TV, KL and
Wasserstein distances.  For instance, by using high-accuracy SLC
samplers, we obtain TV and KL mixing time guarantees for multi-modal
distributions that scale as $\sqrt d,\polylog(1/\varepsilon))$, and
for sampling from a $\kappa$-SLC distribution, the adaptivty of our
scheme reduces the dependence to logarithmic in $\kappa$.




\subsection{Tweedie-based structure of diffused scores}
\label{SecSteinScore}

We now present the background that underlies our development.
Score-based sampling procedures operate on a sequence of random
variables that are transformed by a simple linear operation with
Gaussian noise.  So as to guide the reader, here we lay out the
induced Tweedie structure.  Given a random vector $\Uvar \in \real^d$
and a pair of positive scalars $a$ and $b$, consider the update
\begin{align}
\label{EqnBasicModel}  
  \Vvar & = a \Uvar + b \Wvar \qquad \mbox{where $\Wvar \sim
    \Normal(0, \IdMat)$ is standard Gaussian.}
\end{align}
This transformation is a form of annealing: the density $\pdensv$ of
$\Vvar$ will be smoother than the density $\pdensu$ of $\Uvar$, since
it is obtained by convolving $\pdensu$ with the Gaussian density.

\paragraph{First-order Tweedie and sampling:}  The
classical Robbins--Tweedie
formula~\citep{robbins1956empirical,Miyasawa1961,Efron11} guarantees
that the \emph{Stein score} function $\nabla \log \pdensv$ can be
expressed in terms of \mbox{$\Exs[\Uvar \mid \Vvar = v]$,} which
allows for denoising-based score estimation.  Knowledge of this score
function allows us to draw samples from $\pdensv$ using 
gradient-based sampling procedures.  For our scheme, it is essential
that the \emph{conditional score} $\nabla_u \log \pucondv$ also has a
simple representation.  In particular, we have
\begin{align}
\nabla_u \log \pucondv(u \mid v) & = \nabla_u \log \pdensu(u) +
\nabla_u \log \pvcondu(v \mid u) \; = \; \nabla_u \log \pdensu(u)
-\frac{a}{b^2} \big(a u - v\big),
\end{align}
where the second equality follows from the fact that $(\Vvar \mid
\Uvar = u) \sim \Normal(a u, b^2 \IdMat)$, and the form of the
Gaussian density.  Consequently, knowledge of the marginal score
$\nabla_u \log \pdensu(u)$ gives us knowledge of the conditional
score.  In particular, we can then use gradient-based algorithms to
draw samples from the backwards conditional distribution $\pucondv$.
In summary, for the $1$-step model $V = a U + b W$, knowledge of the
marginal score functions enables us to exploit fast algorithms for
both (a) generating samples from the marginal distributions $\pdensu$
and $\pdensv$, and (b) generating samples from the backward
conditional $\pucondv$.

\paragraph{Second-order Tweedie and Hessian structure:}
Instead of focusing on the score function---that is, the first
derivative of the log density---the bulk of our analysis is instead
focused on second derivatives.  More precisely, still focusing on the
update $V = a U + b W$, we introduce the two Hessian matrices
\begin{align}
\label{EqnDoubleHessian}  
  \HessU(u) \defn - \nabla^2 \log \pdensu(u) \quad \mbox{and} \quad
  \HessV(v) \defn - \nabla^2 \log \pdensv(v), 
\end{align}
associated with the marginal distributions $\pdensu$ and $\pdensv$
over $\Uvar$ and $\Vvar$ respectively, along with the Hessian
$\BackHessU(u, v) \defn - \nabla_u^2 \log \pdensback(u \mid v)$
associated with the conditional distribution of $\Uvar \mid \Vvar$.
To be clear, our analysis makes central use of these second-order
objects, but the standard sampling schemes that we use to solve
sub-problems are still based on first-order information only.  The
following standard results, proved for completeness in
Appendix~\ref{AppProofLemForwardBackward}, play a central role in our
analysis:
    \begin{subequations}    
      \begin{align}
        \label{EqnForwardHessOne}
\hspace*{-0.5in} \mbox{\bf{Second-order Tweedie:}} & \qquad
\underbrace{\HessV(v)}_{- \nabla^2 \log \pdensv(v)}  = \frac{1}{b^2}
\Big \{ \Id - \frac{a^2}{b^2} \cov( \Uvar \mid \Vvar = v) \Big \},
\qquad \mbox{and} \\
        \label{EqnBackwardHessOne}
        \hspace*{-0.05in} \mbox{\bf{Backward conditional Hessian:}} & \qquad
\quad \underbrace{\BackHessU(u, v)}_{-\nabla^2_u \log \pdensback(u
  \mid v)} = \HessU(u) + \frac{a^2}{b^2} \Id.
      \end{align}
    \end{subequations}


\section{Statement of main results}
We now turn to the statement of our main results, including
exponentially accelerated log-concave sampling (\Cref{ThmConcave}
in~\Cref{SecConcave}), and guarantees for general multi-modal
distributions (\Cref{ThmMultiModal} in~\Cref{SecMulti}).

\subsection{Exponential acceleration  for log-concave sampling}
\label{SecConcave}

We begin by studying the consequences of our modular scheme for the
problem of sampling from a smooth and strongly log-concave (SLC)
distribution.  For a given pair of scalars $0 < \mlower \leq \mupper <
\infty$, we say that a twice differentiable density $\ptarget$ is
$(\mlower, \mupper)$-SLC if its negative log Hessian satisfies the
sandwich relation $\mlower \IdMat \; \preceq - \nabla^2 \log
\ptarget(x) \; \preceq \; \mupper \IdMat$ for all $ x\in \real^d$.
The ratio $\kappa \defn \mupper/\mlower \geq 1$ defines the
\emph{condition number} of the problem.  While standard procedures
exhibit iteration complexity scaling linearly in $\kappa$, the main
result of this section gives a simple modular scheme that, by
exploiting the availability of diffused score functions, provides an
\emph{exponential acceleration}, in particular reducing the dependence
to $\log \kappa$.

\paragraph{Procedure and a general guarantee:}

Given a trajectory length $K$ and initial value $Y_0 = X$, consider
the Gaussian forward path given by
\begin{align}
\label{EqnSequenceConcave}
Y_{k+1} = a_k Y_k + \sqrt{1 - a_k^2} \: W_k \qquad
\text{for $k = 0, 1, \ldots, K - 1$,}
\end{align}
where $W_k \sim \Normal(0, \IdMat)$, and $a_k \in (0,1)$ are stepsizes
to be chosen.  The \emph{sub-problems} to be solved in traversing the
backward path are (i) sampling the terminal distribution $Y_K \sim
p_K$, and (ii) for each \mbox{$k = K-1, \ldots, 0$,} sampling from the
backward conditionals $p_{k \mid k+1}$ of $Y_k \mid Y_{k+1}$.  Letting
$\MyProb_0^K \defn \{p_K \} \cup \{p_{k \mid k+1} \}_{k=0}^{K-1}$
denote the collection of all such sampling sub-problems, the following
result guarantees the existence of a ``short'' trajectory such that
all sub-problems $\MyProb_0^K$ are strongly log-concave with condition
number at most $2$:
\mygraybox{
\begin{theorem}[Logarithmic reduction to SLC black box sampling]
\label{ThmConcave}    
Given any $(\mlower, \mupper)$-strongly log-concave target density
$p$, there is a forward trajectory~\eqref{EqnSequenceConcave} of
length at most
\begin{subequations}
\begin{align}
\label{EqnConcaveRounds}
  K + 1& \leq 2 + \log_2 \big( \mupper/\mlower \big) \; = \; 2 +
  \log_2(\kappa),
\end{align}
such that all sampling sub-problems $\MyProb_0^K$ are SLC with
condition number at most $2$.  Consequently, for each distance
\mbox{$D \in \{\dtv, \dkl, \sqrt{\mlower} \Wass_2 \}$,} given
\emph{any} black-box SLC sampler with query complexity $\NitSLC$, we
can perform \mbox{$\varepsilon$-accurate} sampling (in distance $D$)
from the target $p$ in at most
    \begin{align}
 \label{EqnTotalConcave}
 T(\varepsilon) & = \sum_{k=0}^K \NitSLC \big(\varepsilon/(K+1) \big)
 \; \leq \; \big \{2 + \log_2 (\kappa) \big \} \; \NitSLC
 \Big(\frac{\varepsilon}{2 + \log(\kappa)} \Big) \qquad
 \mbox{queries.}
    \end{align}
\end{subequations}
\end{theorem}
}
In the proof, given in~\Cref{SecConcaveAnal}, we give an explicit but
adaptively chosen sequence $\{a_k\}_{k=1}^K$ that certifies the
theorem's claims. Analyzing the resulting sequence exploits the
second-order Tweedie formula~\eqref{EqnForwardHessOne} as well as the
backward Hessian structure~\eqref{EqnBackwardHessOne}.  The key
technical challenge is showing the logarithmic
scaling~\eqref{EqnConcaveRounds}.  We do so via a combination of the
Cram\'{e}r--Rao bound and the Brascamp--Lieb bound to sandwich the
SLC-conditioning $(m_k, M_k)$ of the intermediate problems along the
sequence~\eqref{EqnSequenceConcave}; we
highlight~\Cref{LemSpectralControl} as a key result.

\subsubsection{Some specific consequences}

\Cref{ThmConcave} is a general reduction that allows for any black-box
SLC sampler to be used in solving the sub-problems.  Specific
consequences of the iteration complexity bound~\eqref{EqnTotalConcave}
are easy to extract for a range of samplers; in all cases, we obtain a
$\polylog(\kappa)$ dependence, whereas the $(d,
\varepsilon)$-dependence changes with the chosen
sampler.\footnote{Since the sub-problems all have condition number at
most $2$, we need only focus on their dependence on the pair $(d,
\varepsilon)$.}

\paragraph{Low-accuracy samplers:}

We begin the low-accuracy samplers, meaning that the iteration
complexity scales with $1/\varepsilon$.  Recall
from~\Cref{SecSteinScore} that access to annealed score functions
gives us access to the required gradients for both marginal and
backward conditional sampling.  Using the ULA updates, it is known
that it suffices to take $\Nula(\varepsilon) \asymp d/\varepsilon^2$
steps to achieve $\varepsilon$-accuracy in KL; combined with our
general guarantee, this yields a total iteration complexity of
\mbox{$T(\varepsilon) \asymp \frac{d\log^3(\kappa)}{\epsilon^2}$.}
The linear in $d$-dimension dependence can be reduced via the use of
faster low-accuracy samplers; for instance, schemes based on
underdamped Langevin would achieve $\sqrt{d}$-scaling, so an overall
complexity \mbox{$T(\varepsilon) \asymp \frac{\sqrt{d} \,
    \log^3(\kappa)}{\epsilon^2}$} in KL distance.  As another example,
if we use randomized midpoint
discretization~\citep{ShenLee2019RandomizedMidpoint}, then we obtain a
total iteration complexity $T(\varepsilon) \asymp
\frac{d^{1/3}}{\varepsilon^{2/3}} \log^{5/3}(\kappa)$ in the
$\Wass_2$-distance.
  
\paragraph{High-accuracy guarantees:}  Turning to high-accuracy samplers,
as previously discussed, there are Metropolized
samplers~\citep{altschuler_chewi_2023,Chewi_book} with KL iteration
complexity $\Nfast(\varepsilon) \asymp \sqrt{d}
\polylog(1/\varepsilon)$. Coupled with the
bound~\eqref{EqnTotalConcave}, this yields an overall iteration
complexity scaling as
\begin{align}
\label{EqnSLCFast}  
T(\varepsilon) & \asymp \sqrt{d} (1 + \log(\kappa)) \log(\frac{1 +
  \log(\kappa)}{\varepsilon}).
\end{align}
To be clear, classical Metropolized samplers require access to both
the score function (first derivative), and log density values; in
certain cases, log densities can be estimated
(e.g.,~\cite{guth2025learning}), but this is not always the case.  In
their work on RTK-MALA, \cite{ZhangNew} tackled this issue, and showed
how to approximate Metropolis sampling using only gradient
information, but did not achieve $\Nfast$ complexity.  In concurrent
work, \citet{Sasha} analyze a ``Bernoulli-factory'' randomized
approximation, and show how it can be leveraged to implement a purely
gradient-based sampler that achieves iteration complexity $\Nfast$.
Since the reduction in~\Cref{ThmConcave} allows \emph{for any SLC
sampler}, we can adopt their particular scheme; doing so leads to a
purely score-based sampler that achieves the
guarantee~\eqref{EqnSLCFast}.  This bootstrapping via our scheme
yields a substantial refinement of their initial guarantee, sharpening
the linear in $\kappa$ scaling from their guarantee to logarithmic in
$\kappa$.

\subsubsection{Robustness to score errors}

Letting $p_k$ denote the marginal distribution of $Y_k$ at round $k$,
we have stated~\Cref{ThmConcave} for samplers that operate using the
exact score functions $\score{k}(x) \defn \nabla_x \log p_k(x)$.  In
practice, these exact score functions are not available, but are
instead estimated using samples (e.g., via denoising using the
first-order Tweedie representation).  Here we describe how our
guarantees remain robust when the procedure is implemented using
estimates $\scorehat{k}$ of the true score functions.  In addition to
the computational error tracked by~\Cref{ThmConcave}, this guarantee
involves additional error in terms of the score differences
$\scorehat{k} - \score{k}$.

For concreteness, we focus on the TV error, and state a general
stability result for score errors, one that can be applied both
to~\Cref{ThmConcave} and~\Cref{ThmMultiModal} to follow in the sequel.
At time $k$, let \mbox{$\Qback(\cdot \mid y)$} denote the family of
backward transition kernels defined by the score estimate
$\scorehat{k}$, and let \mbox{$\Qbacktil_k(\cdot \mid y)$} be the
backwards transition defined by a sampler used to implement this
backward transition.  Suppose that: (a) the sampler is
$\varepsilon_k$-accurate, meaning that \mbox{$\dtv(\Qbacktil_k(\cdot
  \mid y), \Qback_k(\cdot \mid y)) \leq \varepsilon_k$} for each $y$;
and (b) the true backward transition $\Pback_k(\cdot \mid y)$ has
density that is $\scon_k$-strongly log-concave.

\mygraybox{
  \begin{lemma}[Robustness to score errors]
\label{LemRobust}    
Under the stated conditions, for each backward step $k = 0, 1, \ldots,
K-1$, the marginal $\qtil_k$ satisfies the TV norm recursion
  \begin{align}
\label{EqnTVRobustScore}    
\dtv(\qtil_k,p_k) \leq \varepsilon_k + \frac{\|\scorehat{k} -
  \score{k} \|_{L^2(\qbar_k)}}{2\sqrt{\scon_k}} + \dtv(\qtil_{k+1},
p_{k+1}) \quad \mbox{where $\qbar_k = \qtil_{k+1} \Qback_k$.}
  \end{align}
\end{lemma}
}
\noindent See~\Cref{SecProofLemRobust} for the proof of this claim. \\

The TV guarantee in~\Cref{ThmConcave} makes use of the
bound~\eqref{EqnTVRobustScore} with no score error ($\scorehat{\ell} =
\score{\ell}$, and the settings $\varepsilon_\ell = \varepsilon/(K +
1)$, so that the total error scales as $\sum_{\ell=0}^K
\varepsilon_\ell = \varepsilon$.  Since~\Cref{LemRobust} allows for
score error,\footnote{The lemma applies to the backward transitions,
and we can also prove an analogous stability result for the terminal
sub-problem $p_K$.} we can prove a generalized result that involves
the total error
\begin{align}
 \label{EqnTotalError}
\sum_{\ell=0}^K \varepsilon_\ell + \sum_{\ell=0}^K
\frac{\|\scorehat{\ell} - \score{\ell}
  \|_{L^2(\qbar_\ell)}}{2\sqrt{\scon_\ell}},
\end{align}
where, in the setting of~\Cref{ThmConcave}, the strong log-concavity
parameter is given by $\scon_\ell = 2 + 2^{-\ell} (\kappa - 1)$, where
$\kappa > 1$ is the condition number.



\subsection{Multi-modal setting}
\label{SecMulti}

We now turn to sampling from a general multi-modal distribution.  In
this case, given some $\sigtar > 0$ that is user-specified, our goal
is to draw samples from the random vector $Z \defn \Xdata + \sigtar
W_0$, where $W_0 \sim \Normal(0, \Id)$. As in diffusion analyses, we
refer to the parameter $\sigtar$ as the early stopping error.  With
this parameter fixed, our goal is to develop algorithms that produce
$\varepsilon$-accurate samples from the distribution $p_Z$ of $Z$.  In
order to do so, it is convenient to work with the rescaled random
vector $Y_1 = \frac{Z}{\sqrt{2} \, \sigtar}$.  This rescaling has no
effect on the distances $\dkl$ and $\dtv$.

\subsubsection{Adaptive scheme and its analysis}
Having reduced our problem to sampling from $p_1 \equiv p_{Y_1}$, we
note that $Y_1$ can be written as $Y_1 \defn \tfrac{1}{\sqrt{2}} X +
\tfrac{1}{\sqrt{2}} W_0$, \mbox{where $X \defn \Xdata/\sigtar$} is a
rescaled version of the original variable $\Xdata$.  In our analysis,
we state conditions directly on $X$.  Given this $Y_1$, we then
generate the forward sequence
\begin{subequations}
\begin{align}
  \label{EqnYrecursion}
  Y_{k+1} = a_k Y_k + \sqrt{1 - a_k^2} W_k, \qquad \text{for stepsizes
    $a_k \in (0,1)$ to be chosen adaptively.}
\end{align}
Our adaptive choice of stepsizes depends on the sequence
\begin{align}
\label{EqnDefnCovCon}
\CovCon_k \: \defn \: \sup_{y \in \real^d} \opnorm{\cov(X \mid Y_k = y)},
\end{align}
defined for each $k = 1, 2, \ldots$ in the forward trajectory.  We
construct the forward sequence with a stepsize sequence $\{a_k\}_{k
  \geq 0}$ based on the initialization $a_0 = 1/\sqrt{2}$, and for $k
= 1, 2, \ldots$, the adaptive updates
\begin{align}
\label{EqnMultiStepsize}
\IterB_k \: \defn \: 4 \, \CovCon_k \prod_{\ell = 0}^{k-1} a_\ell^2 \quad
\mbox{and} \quad a_k^2  = \frac{2 \IterB_k + 2}{2 \IterB_k + 3}.
\end{align}
\end{subequations}
As before, for a trajectory length $K$, we let $\MyProb_0^K \defn
\{p_K \} \cup \{p_{k \mid k+1} \}_{k=0}^{K-1}$ be the collection of
all sampling sub-problems: the terminal distribution $p_K$, and each
of the backward conditional distributions.

\mygraybox{
\begin{theorem}[SLC reduction for multi-modal case]
\label{ThmMultiModal}
Given the adaptive stepsize choice~\eqref{EqnMultiStepsize}, consider
any trajectory $K$ length such that $\prod_{\ell = 0}^{K -1} a_\ell^2
\leq \frac{1}{8 \CovCon_K}$.  Then we are guaranteed that all of the
sampling sub-problems $\MyProb_0^K$ are strongly log-concave (SLC)
with condition number at most $4$.  Consequently, for any distance $D
\in \{\dtv, \dkl \}$, by using \emph{any} SLC sampler with query
complexity $\NitSLC$, we can obtain $\varepsilon$-accurate samples
using
\begin{align}
\label{EqnMultiFinal}  
  T(\varepsilon) = \sum_{k=0}^{K} \NitSLC \big( \varepsilon/(K+1)
  \big) \qquad \text{ queries.}
\end{align}
\end{theorem}
}
\noindent See~\Cref{SecProofThmMultiModal} for the proof.  Here the
main technical challenge is showing how the adaptive stepsize
choice~\eqref{EqnMultiStepsize} suffices to ensure the SLC property of
the full collection $\MyProb_0^K$ of sub-problems.  As
with~\Cref{ThmConcave}, this reduction allows for any SLC sampler, and
so has consequences for both low- and high-accuracy schemes.

\paragraph{Low-accuracy guarantees:}  If we use ULA as our base
sampler, then we achieve an overall guarantee of $T(\varepsilon)
\asymp K^3 \frac{d}{\varepsilon^2}$; this is comparable to the RTK-ULA
guarantee in the paper~\citep{ZhangNew}.  On the other hand, we could
also use underdamped Langevin to achieve the scaling $T(\epsilon)
\asymp K^3 \frac{\sqrt{d}}{\varepsilon^2}$. If the trajectory length
$K$ is independent of dimension, then the end-to-end procedure
exhibits $\sqrt{d}$-dependence, as opposed to the optimal linear
$d$-scaling that arise from optimal KL bounds for ULA-based
discretizations of
diffusions~\citep{conforti2023score,benton2024nearly}.  We can also
use a more sophisticated Langevin discretization to go below the
$\sqrt{d}$ barrier; for instance, if we use a randomized midpoint
scheme (RMC), then known TV-guarantees for
RMC-Langevin~\citep{gupta_cai_chen_2025} give us an overall complexity
of $T(\varepsilon) \asymp K^{7/3} d^{5/12}/\varepsilon^{4/3}$.

\paragraph{High-accuracy results:}
As discussed following~\Cref{ThmConcave}, we can achieve
$\varepsilon$-accurate in KL samples for each sub-problem with
complexity $\Nfast(\varepsilon) \asymp \sqrt{d}
\polylog(1/\varepsilon)$, either by using Metropolized
schemes~\citep{altschuler_chewi_2023} that use gradients and density
values, or by applying the gradient-only SLC-Metropolis
sampler~\citep{Sasha}, as discussed following~\Cref{ThmConcave}. In
either case, when combined with the guarantee~\eqref{EqnMultiFinal},
our reduction yields an overall iteration complexity of
\begin{align}
\label{EqnMultiBoundKL}  
  T(\varepsilon) & = \Otil \Big( K \, \sqrt{d} \, \polylog
  \big(\frac{K}{\varepsilon} \big) \Big).
\end{align}
In general, the trajectory length $K$ is a function of the full
sequence $\CovCon_k$ from equation~\eqref{EqnDefnCovCon}; in the next
section, we bound $K$ under a worst-case assumption.

\subsubsection{Trajectory length under worst-case assumptions}
\label{SecMultWorst}

We now turn to an analysis of the trajectory length $K$ under a
particular worst-case assumption, namely that the only knowledge about
the sequence $\{\CovCon_k \}_{k \geq 1}$ is the existence of some
$\Cmax < \infty$ such that
\begin{align}
  \label{EqnWorstCase}
\CovCon_k \: \defn \: \sup_{y \in \real^d} \opnorm{\cov(X \mid Y_k =
  y)} \leq \Cmax \qquad \mbox{for all $k = 1, 2, \ldots$,}
\end{align}
but no further structure is given.  By comparison to the second-order
Tweedie formula~\eqref{EqnForwardHessOne}, we see that this type of
uniform bound corresponds to assuming the score function has a bounded
Lipschitz constant.  By careful analysis of the adaptive
sequence~\eqref{EqnMultiStepsize}, we obtain the following guarantee:
\mygraybox{
  \begin{corollary}
    \label{CorMultiWorst}
Under the uniform bound~\eqref{EqnWorstCase}, the guarantees
of~\Cref{ThmMultiModal} apply with trajectory length at most $\Kfinal
\leq 7 (1 + \Cmax)$.
  \end{corollary}
}
\noindent See Appendix~\ref{SecProofCorMultiWorst} for the proof.

Substituting the upper bound on $K$ from Corollary~\ref{CorMultiWorst}
into our earlier bound~\eqref{EqnMultiBoundKL}, we find a worst-case
iteration complexity scaling as
\begin{align}
\label{EqnTworst}
\Tworst(\varepsilon) & \asymp \Cmax \sqrt{d} \polylog(1/\varepsilon),
\end{align}
where $\Cmax$ is the Lipschitz constant.  Comparing to past
results~\citep{chen2023improved,chen2023sampling} on discretized
diffusions, these bounds exhibit a quadratic dependence on a Lipschitz
condition, and $d/\varepsilon^2$ dependence on $(d, \varepsilon)$, as
opposed to $\Cmax \sqrt{d} \polylog(1/\varepsilon)$ in the
guarantee~\eqref{EqnTworst}.  The $\polylog(1/\varepsilon)$-scaling
for diffusion models first appeared in RTK-MALA guarantee
of~\cite{ZhangNew}, albeit with a substantially worse $d^2$-scaling,
and applicable only to TV distance.  In concurrent work, \citet{Sasha}
impose a related Lipschitz assumption, and obtain guarantees with the
same $\sqrt{d} \: \polylog(1/\varepsilon)$ scaling using a standard
diffusion approach. Their result imposes high-probability Lipschitz
control as opposed to the worst-case bound imposed in
Corollary~\ref{CorMultiWorst}.  It would be interesting to use our
framework to obtain bounds on $K$ depending on a milder geometric
quantity; we strongly suspect this should be possible, since the
sequence $\{\CovCon_k\}_{k \geq 1}$ evolves in a very structured way.


\subsubsection{Some open questions and extensions}
\label{SecOpen}

Our scheme and results suggest a variety of other open questions and
extensions.  Let us comment on a few of them here.

\paragraph{Best of all worlds?}   It is interesting to make
further explicit comparisons in a specific setting, which we refer to
as the \emph{bounded $(R, \sigma)$-model.}  Suppose that we wish to
sample $Y_1 = \Xtil + \sigma W$, where $\|\Xtil \|_2 \leq R$.  In this
case, the worst-case bound~\eqref{EqnWorstCase} holds with $\Cmax =
R^2/\sigma^2$, and our approach yields the overall iteration
complexity $T_{R, \sigma}(\varepsilon) \asymp (R/\sigma)^2 \sqrt{d}
\polylog(1/\varepsilon)$.  In work prior to our own, state-of-the-art
results on discretized
diffusions~\citep{conforti2023score,benton2024nearly,li_wei_chi_chen_2024}
provide iteration complexities with $d/\varepsilon^2$ scaling,
inferior to the $\sqrt{d} \log(1/\varepsilon)$ scaling in $T_{R,
  \sigma}(\varepsilon)$, but with a milder poly-logarithmic dependence
on $R/\sigma$.  We are thus led to a natural open question: is it
possible to obtain a ``best-of-both-worlds'' guarantee: i.e., with
poly-logarithmic dependence on both $(R/\sigma)$ and $(1/\epsilon)$?

\paragraph{Reductions beyond the SLC class:}  \Cref{ThmConcave} and~\Cref{ThmMultiModal} both exploit
trajectories of strongly log-concave (SLC) problems, which allows us
to exploit fast SLC sampling algorithms.  However, there are many
other classes of distributions for which fast samplers are available,
including those satisfying geometric relations such as log-Sobolev
(LSI) or Poincar\"{e} inequalities. For instance, the class of LSI
distributions is a strict superset of the SLC class, allowing for
considerable multi-modality depending on the LSI constant.  There are
various fast samplers available for LSI distributions
(e.g.,~\cite{vempala_wibisono_2022,CheGam23,altschuler_chewi_2023}),
so that the scope of our scheme could be enlarged substantially by
instead reducing to a sequence of well-conditioned LSI problems.  One
natural open question is whether (for $1$-Lipschitz problems) one can
obtain sampling guarantees with logarithmic dependence on the LSI
constant, thereby generalizing the logarithmic dependence on condition
number $\kappa$ in~\Cref{ThmConcave}.

\comment{
\paragraph{Milder dimension dependence:}  An attractive feature of the general
reduction of~\Cref{ThmMultiModal} is that it allows for the use of any
sampler for the SLC problems that define the backward trajectory.  In
deriving concrete results thus far, we have assumed the use of
sampling routines with $\sqrt{d}$-dependence, so that our summary
results include this scaling.  However, it is worth noting that there
are various first-order schemes with even milder dimension dependence
when additional structure is available. For instance, under a
higher-order smoothness condition, Mou et al.~\cite{MouHighOrder21}
exhibit a sampling algorithm, one that uses only first-order
information but makes use of a lifted form of dynamics, with $d^{1/4}$
scaling. On the other hand, Chen and Gamirty~\cite{CheGam23} provide a
refined analysis of the MALA algorithm that yields square-root
dependence on the Hessian trace, as opposed to $\sqrt{d}$.  Any of
these improved samplers can be used within our modular scheme, and it
inherits any improved dimension dependence in the final iteration
complexity.

\paragraph{Manifold and other low-dimensional structures:}
As noted in~\Cref{SecRelated}, there is a rapidly evolving line of
work on diffusion sampling when the original distribution exhibits
some type of low-dimensional structure
(e.g.,~\cite{debortoli2022convergence,pidstrigach2022score,debortoli2022riemannian,chen2023score,azangulov2024convergence}),
with recent results giving linear scaling in effective
dimension~\cite{potaptchik2024linear,yuting2024_lowdim}. It would be
interesting to prove analogous results for our modular scheme, so as
to benefit both from the poly-logarithmic scaling in
$(1/\varepsilon)$, and adaptivity to low-dimensional structure.  For
our scheme, one key question to address is the iteration complexity of
samplers for SLC distributions obtained by annealed smoothing of a
low-dimensional distribution.
}



\section{Proof sketches}

We now provide sketches of our two main results: ~\Cref{ThmConcave}
in~\Cref{SecConcaveAnal} and~\Cref{ThmMultiModal}
in~\Cref{SecProofThmMultiModal}.

\subsection{Proof sketch of~\Cref{ThmConcave}}
\label{SecConcaveAnal}
At a high level, the proof consists of three main steps:
\begin{enumerate}[leftmargin=*, itemsep=0pt]
\item[(1)] \underline{Rescaling:} We apply a rescaling argument to
  reduce the problem to a sampling problem that $(1,
  \mupper/\mlower)$-strongly log-concave; see
  Appendix~\ref{AppRescaling} for this argument.
\item[(2)] \underline{Adaptive stepsizes:} We prescribe an adaptive
  choice of stepsizes $\{a_k\}_{k=0}^{K-1}$, and analyze the evolution
  of marginal distributions $p_k \equiv p_{\myspecial{y_k}}$ of the
  forward process, and the conditional distributions $p_{k \mid k +1}
  \equiv p_{\myspecial{y_k} \mid \myspecial{y_{k+1}}}$ of the backward
  process; and
\item[(3)] \underline{Stability analysis:} We control error
  propagation throughout the entire backward process;
  see~\Cref{SecConcaveError} for details of this step.
\end{enumerate}
In the main body, we provide more detail on adaptive stepsize and
trajectory analysis in Step 2, and the stability analysis in Step 3.

\subsubsection{Forward recursion with adaptive stepsize choice}
\label{SecConcaveAdaptive}
Introduce the shorthand $\mupper_0 = \mupper/\mlower$ for the
smoothness constant of the target distribution after the rescaling
step.  We run the forward recursion~\eqref{EqnSequenceConcave} with
the following choice of stepsizes $\{a_k \}_{k \geq 1}$. Given the
initial value $\mupper_0$, we initialize with $a_0^2 = \mupper_0/(1 +
\mupper_0)$, and then evolve the pair $(\mupper_k, a_k)$ according to
the recursion
\begin{align}
  \label{EqnBasicRecursion}  
  \mupper_{k+1} & \stackrel{(i)}{=} \frac{\mupper_k}{\mupper_k (1 -
    a_k^2) + a_k^2}, \quad \mbox{and} \quad a_k^2 \stackrel{(ii)}{=}
  \frac{\mupper_k}{1 + \mupper_k} \qquad \mbox{at each round.}
\end{align}
Given these stepsize choices, the key technical steps in our proof are
to establish that after $K \leq 1 + \log_2 \mupper_0$ rounds, the
Hessian $\Hess{K}$ satisfies
\begin{subequations}
\begin{align}
  \label{EqnForwardControlOne}
  \IdMat \; \preceq \; \Hess{K}(y) \; \preceq \; 2 \IdMat \qquad
  \mbox{uniformly in $y \in \real^d$,}
\end{align}
and at each round $k \in [K-1]$, the Hessian $\BackHess_k$ of the
backward conditional $p_{k \mid k+1}$ satisfies the sandwich
\begin{align}
  \label{EqnBackwardControlOne}
  1 + \mupper_k \; \preceq \; \BackHess_k(y) \; \preceq \; 2
  \mupper_k \qquad \mbox{uniformly in $y \in \real^d$},
\end{align}
\end{subequations}
Note that the sandwich relation~\eqref{EqnForwardControlOne} ensures
that the terminal distribution $p_K$ is SLC with condition number at
most $2$, whereas the sandwich relations~\eqref{EqnBackwardControlOne}
guarantee that each of the backward conditionals $p_{k \mid k+1}$ are
SLC with condition number at most $2$.  Together, these results
certify that all the sampling sub-problems $\MyProb_0^K$ satisfy the
claims of~\Cref{ThmConcave}.

\subsubsection{Propagation of spectral control}

Proving the claims~\eqref{EqnForwardControlOne}
and~\eqref{EqnBackwardControlOne} hinges on the following auxiliary
result.  Each step in the sequence~\eqref{EqnSequenceConcave} is of
the generic form $\Vvar = a \Uvar + b \Wvar$, where $\Wvar \sim
\Normal(0, \IdMat)$ is standard Gaussian.  Recalling the Hessian
matrices from equation~\eqref{EqnDoubleHessian}, the following result
shows how control of the spectrum of $\HessU(u)$ yields spectral
control of $\HessV(v)$.
\mygraybox{
  \begin{lemma}[Propagation of spectral control]$\quad$
    \label{LemSpectralControl}
    \begin{subequations}
    \begin{enumerate}[leftmargin=*, itemsep=0pt, nosep]
    \item[(a)] Suppose there is some $\mlower \geq 0$ such that
      $\HessU(u) \succeq \mlower \IdMat$ uniformly in \mbox{$u \in
      \real^d$}. Then we have
      \begin{align}
        \label{EqnVlow}
        \HessV(v) & \succeq \; \frac{\mlowerbar}{\mlowerbar b^2 + a^2}
        \Id \qquad \mbox{uniformly in $v \in \real^d$.}
      \end{align}
    \item[(b)] Suppose there is some $\mupper < \infty$ such that
      $\HessU(u) \preceq \mupper \IdMat$ uniformly in $u \in
      \real^d$. Then we have
      \begin{align}
        \label{EqnVup}
        \HessV(v) & \preceq \; \frac{\mupper}{\mupper b^2 + a^2} \Id
        \qquad \mbox{uniformly in $v \in \real^d$.}
      \end{align}      
    \end{enumerate}
    \end{subequations}
\end{lemma}
}
\noindent
See Appendix~\ref{SecProofLemSpectralControl} for the proof. Part (a)
requires strong log-concavity (SLC), and exploits the Brascamp--Lieb
inequality applicable in this case.  Part (b) applies in general, and
makes use of the Cram\'{e}r-Rao lower bound on location estimation.

\subsubsection{Backwards error propagation}
\label{SecConcaveError}

Finally, we provide a lemma that allows us to control propagation of
errors in moving along the backward path.  On one hand, we have the
Markov kernel
\begin{subequations}
\begin{align}
\label{EqnDefnBackMark}
\Pback_k(r)(\cdot) & \defn \int_{\real^d} \pkback(\cdot \mid y)
r(y) d y
\end{align}
that defines the backwards evolution of the true marginals $p_{k+1}
\rightarrow p_k$.  Our backwards sampler defines a second kernel
$\Qback_k$ that underlies the backwards evolution $q_{k+1} \rightarrow
q_k$ of the algorithm's marginals.  By our set-up, with a sufficient
number of iterations, we can assume that our backwards sampler is
$\delta_k$-accurate uniformly in its inputs, meaning that
\begin{align}
  \label{EqnSamplerAccuracy}
D( \Qback_k(e_y), \Pback_k(e_y)) \leq \delta_k \qquad \text{for all $y
  \in \real^d$,}
\end{align}
\end{subequations}
\mygraybox{
  \begin{lemma}[Error propagation in backwards kernel $\Pback_k$]
\label{LemConcaveError}
Under the condition~\eqref{EqnSamplerAccuracy}, for any of the
distances $D \in \{\dtv, \dkl, \Wass_2 \}$, we have
\begin{align}
  \label{EqnSLCStable}
  \underbrace{D(q_k, p_k)}_{\mbox{Round-$k$ error}} & \leq \quad
  \delta_k \quad + \underbrace{D(q_{k+1},
    p_{k+1}).}_{\mbox{Round-$(k+1)$ error}}
\end{align}
\end{lemma}
}

\noindent See Appendix~\ref{SecProofLemConcaveError} for the proof.
For the TV distance $\dtv$, the bound~\eqref{EqnSLCStable} is a
special case of~\Cref{LemRobust} without score error.  Proof of the
bound~\eqref{EqnSLCStable} for the KL divergence $\dkl$ is an easy
consequence of the data processing inequality; the corresponding proof
for Wasserstein is more involved (cf.~\Cref{LemWassStable} for
details).

\subsubsection{Combining the pieces}
\label{SecConcavePieces}

Let us now combine the pieces so as to complete the proof
of~\Cref{ThmConcave}.  Recall that we say that a problem is $2$-SLC if
it is strongly log concave (SLC) with condition number at most $2$.
Let $(p_K, p_{K-1}, \ldots, p_1, p_0)$ denote the true sequence of
marginal distributions, and let $(q_K, q_{K-1}, \ldots, q_1, q_0)$
denote the sequence of distributions generated by our approximate
sampling algorithm. We assume that:
\begin{itemize}[leftmargin=*, itemsep=0pt]
\item At the terminal stage $K$, we generate samples from a
  distribution $q_K$ that is $\varepsilon/K$-close to $p_K$.  From
  the guarantee~\eqref{EqnForwardControlOne}, this terminal stage
  problem is $2$-SLC-controlled, so that doing so using our black-box
  SLC procedure requires $\NitSLC(\varepsilon/K)$ calls.
\item In moving backward from $(k+1)$ to $k$, we run enough iterations
  of the sampler so that the bound~\eqref{EqnSamplerAccuracy} holds
  with $\delta_k = \varepsilon/(K+1)$.  From the
  guarantee~\eqref{EqnBackwardControlOne}, each backward conditional
  is $2$-SLC-controlled, so that doing so requires $\NitSLC(
  \varepsilon/(K+1))$ calls.
\end{itemize}
By construction, the total number of calls required is $\sum_{k=0}^K
\NitSLC(\varepsilon/(K+1))$, matching the
claim~\eqref{EqnTotalConcave}.  For a given distance $D \in \{\dtv,
\kull, \Wass_2 \}$, let us compute the error $D(q_0, p_0)$ at the
initial stage.  By expanding the bound~\eqref{EqnSLCStable} at each
round with $\delta_k = s_k \varepsilon$, we have $D(q_0, p_0) \; \leq
\; \varepsilon/(K+1) + D(q_{1}, p_1) \leq \frac{1}{K} \sum_{k=1}^K
\varepsilon = \varepsilon$, as claimed.

\subsection{Proofs of~\Cref{ThmMultiModal} and Corollary~\ref{CorMultiWorst}}
\label{SecProofThmMultiModal}

We now turn to the proof sketch for~\Cref{ThmMultiModal}, along with
related Corollary~\ref{CorMultiWorst}. It suffices to prove the result
for $\dkl$, since upper bounds on $\dtvsquared$ follow as a
consequence of Pinsker's inequality. This proof does not require
rescaling step, and we can re-use the backward error propagation for
$\dkl$ given in~\Cref{LemConcaveError}.  (However, the Wasserstein
error propagation there does not apply to the multi-modal case, as the
Wasserstein proof in~\Cref{LemConcaveError} exploits the SLC
structure.)  Consequently, the new technical effort required is
analyzing the evolution of the forward and backward processes induced
by the adaptive stepsize sequence~\eqref{EqnMultiStepsize}.

Recall the Hessian matrices from equation~\eqref{EqnDoubleHessian}.
The following result shows how control of the spectrum of $\HessU(u)$
yields spectral control of $\HessV(v)$.  The following lemma controls
the Hessian structure of the trajectory, both in terms of the marginal
Hessians of the forward process, and the conditional Hessians of the
backward trajectory.  It allows us to show that all relevant sampling
sub-problems satisfy the requisite SLC problems.

\mygraybox{
\begin{lemma}[Trajectory control for multi-modal case]
  \label{LemMultiForward}
  Given the adaptive stepsize choice~\eqref{EqnMultiStepsize},
  the following properties hold:
\begin{enumerate}[leftmargin=*, itemsep=0pt, nosep]
\begin{subequations}    
\item[(a)] At each round $k \in [K]$, we have
\begin{align}
\label{EqnFhessian}  
  - \IterB_k \IdMat \; \preceq (1 - \IterB_k) \IdMat \; \preceq \;
  \PlainHess_k(y) \; \preceq \; 2 \IdMat.
\end{align}    
\item[(b)] For each $k \in[K-1]$, the backwards conditional
  distribution $p_{k \mid k+1}$ has a Hessian $\BackHess_k(y_k)$ that
  satisfies the sandwich
  \begin{align}
\label{EqnBhessian}    
   ( \IterB_k + 2 ) \, \IdMat \;\preceq\; \BackHess_k(y) \;\preceq\; 2 \,
\big( \IterB_k + 2 \big) \, \IdMat \qquad \mbox{for all $y \in \real^d$.}
  \end{align}
\end{subequations}
\end{enumerate}
  \end{lemma}
}
\noindent See Section~\ref{SecProofLemMultiForward} for the proof of
this claim. \\


\paragraph{Completing the proof:}
We now sketch how to complete the proof of~\Cref{ThmMultiModal}.  From
the theorem set-up, recall that the trajectory length $K$ is chosen to
ensure that $\prod_{k=0}^{K -1} a_k^2 \leq \frac{1}{8 \CovCon_K}$.
From the definition of the sequence $\IterB_k$, this bound implies
that $\IterB_K \leq 1/2$.  Thus, from part (a)
of~\Cref{LemMultiForward}, we see that the terminal Hessian $\Hess{K}$
satisfies the sandwich bound $\tfrac{1}{2} \IdMat \preceq \Hess{K}(y)
\preceq 2 \IdMat$ \text{uniformly in $y \in \real^d$,} so that it is
SLC with condition number at most $4$, as claimed.  From part (b)
of~\Cref{LemMultiForward}, we see that the backward Hessians
$\BackHess_k$ have condition number at most $2$.  Thus, we have
established that the adaptive stepsize sequence constructs a sequence
of $K$ sub-problems that are each SLC with condition number at most
$4$.  Finally, the sampling guarantees~\eqref{EqnMultiFinal} given
in~\Cref{ThmMultiModal} follows from the same proof template used
in~\Cref{SecConcavePieces}.


\section{Discussion}
\label{SecDiscussion}

In this paper, we have described a modular approach to sampling from a
given target distribution $\ptarget$ based on the availability of
annealed score functions.  It can be understood as a form of
``divide-and-conquer'', showing how the original sampling problem is
reducible to a $K$-length sequence of sub-problems, each defined by a
well-conditioned and strongly log-concave (SLC) distribution.  Using
this reduction, we proved novel results both for uni-modal and
multi-modal distributions.  For sampling from a SLC distribution
(\Cref{ThmConcave}), we exhibited efficient methods scaling as $\log
\kappa$ with the condition number $\kappa$, and dimension dependence
ranging in $\{d^{1/3}, \sqrt{d}, d \}$, depending on the SLC sampler
used to solve the sub-problems.  Using high-accuracy SLC samplers
leads to $\sqrt{d} \polylog(1/\varepsilon)$ scaling.  We established
similar guarantees for general multi-modal distributions in
\Cref{ThmMultiModal}; here the $\sqrt{d} \polylog(1/\varepsilon)$
scaling obtained by high-accuracy samplers improves upon the best
known results for diffusion samplers prior to our work.

Our work leaves open various questions and potential extensions.  In
Corollary~\ref{CorMultiWorst}, we provided a worst-case bound on the
trajectory length $K$ for multi-modal problems; this crude analysis
does not exploit any structure of the process, and we suspect that it
could be sharpened.  At a broader level, we focused on adaptive
schemes that generate sequences of SLC sub-problems with constant
conditioning, but our general idea is not limited to this choice.  It
could be enriched by instead reducing to richer classes of
distributions for which fast schemes are available (e.g., those
satisfying a log-Sobolev inequality (LSI)).  A related open problem is
whether such a reduction could be used to prove the natural
generalization of~\Cref{ThmConcave}: an efficient sampling scheme with
logarithmic dependence on the LSI constant. \footnote{Assume that the
log distribution is $1$-smooth, so that LSI constant plays the role of
condition number $\kappa$.}  The use of LSI sub-problems could also
yield improvements for multi-modal sampling; in our
current~\Cref{ThmMultiModal}, being forced to ensure SLC on the
backward path can be restrictive.  We also raised a number of open
question associated with the bounded $(R, \sigma)$-model;
see~\Cref{SecOpen} for further discussion.

\acks{ This work was partially supported by a Guggenheim Fellowship and
  grant NSF DMS-2311072 from the National Science Foundation.}

\newpage

{\small{
\bibliography{mjwain_super,martin_papers,merged}
}}




\appendix

\crefalias{section}{appendix} 

\section{Auxiliary results for~\Cref{ThmConcave}}

In this appendix, we collect together the proofs of various
auxiliary results related to~\Cref{ThmConcave}.

\subsection{Rescaling argument}
\label{AppRescaling}

It is convenient, as a first step, to reduce the problem of sampling
from an $(\mlower, \mupper)$-conditioned distribution to an equivalent
one that is \mbox{$(1, \mupper/\mlower)$-conditioned.}  In order to do
so, letting $X \sim p$ be the original distribution, we define the
rescaled vector $Y = \sqrt{\mlower} X$.  If the original distribution
$p$ is $(\mlower, \mupper)$-conditioned, then the rescaled vector $Y$
has a distribution that is $(1, \mupper/\mlower)$-conditioned.

Now suppose that we can generate samples of random vector $\Ytil$
whose distribution is $\varepsilon$-close to that of $Y$ in a given
distance $D$.  We then define \mbox{$\Xtil = \Ytil/\sqrt{\mlower}$,}
and consider the quality of the $\Xtil$ samples as approximations to
the original $X$.  The argument is slightly different, depending on
whether $D$ is the KL distance or the Wasserstein distance.

When $D$ is the KL distance, it is invariant to this linear
transformation, so that we have
\begin{align*}
\kull(p_{\Xtil}, p_{X}) & = \kull(p_{\Ytil}, p_{Y}) \; = \; \varepsilon.
\end{align*}
Consequently, the samples $\Xtil$ are $\varepsilon$-close to $X$ in KL
distance.

When $D$ is the Wasserstein distance, we use $\Wass_2(X, Y)$ denote
the Wasserstein distance between $p_X$ and $p_Y$. In this case, we
have
\begin{align*}
  \sqrt{\mlower} \Wass_2(\Xtil, X) & = \Wass_2 \big( \sqrt{\mlower}
  \Xtil, \sqrt{\mlower} X \big) \; = \; \Wass_2(\Ytil, Y ) \; = \;
  \varepsilon
\end{align*}
so that the samples are $\varepsilon$-close in the rescaled Wasserstein
norm $\sqrt{\mlower} \Wass_2$.  (Note that the theorem statement
involves this rescaled Wasserstein distance.)

Thus, for the remainder of the proof, we study the rescaled problem
with initial values $\mlower_0 = 1$ and $\mupper_0 =
\frac{\mupper}{\mlower}$.  We prove bounds in terms of $\mupper_0$,
and then recall this transformation.

\subsection{Analysis of forward-backward trajectory}

The following lemma summarizes some key properties of the
forward/backward trajectory:

\mygraybox{
  \begin{lemma}[Properties of forward/backward trajectory]
\label{LemSLCRecursion}    
    \begin{enumerate}[leftmargin=*, itemsep=0pt]
      \begin{subequations}
      \item[(a)] We have the Hessian sandwich
        \begin{align}
          \label{EqnSLCSand}
          \IdMat \; \preceq \; \Hess{k}(y) \; \preceq \; \mupper_k
          \IdMat \qquad \mbox{at each round $k \in [K]$.}
        \end{align}
      \item[(b)] After $K \leq 1 + \log_2 \mupper_0$ rounds, the
        Hessian $\Hess{K}$ satisfies
        \begin{align}
          \label{EqnForwardControl}
          \IdMat \; \preceq \; \Hess{K}(y) \; \preceq \; 2 \IdMat
          \qquad \mbox{uniformly in $y \in \real^d$,}
        \end{align}
so that $\sup \limits_{y \in \real^d} \mycond(\Hess{K}(y)) \leq 2$.        
      \item[(c)] At each round $k \in [K-1]$, the Hessian
        $\BackHess_k$ of the backwards conditional $p_{k \mid k+1}$
        satisfies the sandwich
        \begin{align}
          \label{EqnBackwardControl}
          (1 + \mupper_k) \IdMat \; \preceq \; \BackHess_k(y) \;
          \preceq \; 2 \mupper_k \IdMat \qquad \mbox{uniformly in $y
            \in \real^d$},
        \end{align}
so that $\sup \limits_{y \in \real^d} \mycond(\BackHess_k(y)) \leq 2$.
      \end{subequations}
    \end{enumerate}
  \end{lemma}
}

\begin{proof}
  We first prove the Hessian sandwich~\eqref{EqnSLCSand}, in
  particular via induction on the iteration number $k$.  Beginning
  with the base case $k = 0$, the claim holds because the original
  problem is $(1, \mupper_0)$-conditioned by construction.  Suppose
  that the sandwich~\eqref{EqnSLCSand} holds at step $k$, and let us
  prove that it holds at step $k + 1$.

Beginning with the lower bound, we apply the bound~\eqref{EqnVlow}
from~\Cref{LemSpectralControl} with $\mlower = 1$, $U = Y_k$ and
$V = Y_{k+1}$, and $a^2 = a_k^2$ and $b^2 = 1-a_k^2$.  Doing so yields
\begin{align*}
\Hess{k+1}(y_{k+1}) & \succeq \frac{1}{a_k^2 + (1-a_k^2)} \IdMat \; =
\; \IdMat
\end{align*}
as required.  Similarly, we apply the upper bound~\eqref{EqnVup}
with $\mupper = \mupper_k$ and the same choices as above.  Doing so yields
\begin{align*}
\Hess{k+1}(y_{k+1}) & \preceq \frac{\mupper_k}{\mupper_k (1-a_k^2) +
  a_k^2} \IdMat \; = \; \mupper_{k+1} \IdMat,
\end{align*}
using the definition~\eqref{EqnBasicRecursion} of the $(\mupper_k,
a_k)$ recursion.

Now we establish the claim~\eqref{EqnForwardControl}.  Based on the
Hessian sandwich~\eqref{EqnSLCSand} from part (a), the proof of this
claim amounts to showing that $\mupper_K \leq 2$ after at most $K \leq
1 + \log_2 \mupper_0$ rounds.  We claim that the sequence $\mupper_k$
evolves in a very simple way---namely as
\begin{align}
  \label{EqnExpFast}  
  \mupper_{k+1} & = 1 + \tfrac{1}{2} \big( \mupper_k - 1 \big).
\end{align}
To prove this claim, note that since $1 - a_k^2 = \frac{1}{1 +
  \mupper_k}$, the denominator in the
recursion~\eqref{EqnBasicRecursion} is given by
\begin{align*}
  \mupper_k (1 - a_k^2) + a_k^2 = \frac{\mupper_k}{1 + \mupper_k} +
  \frac{\mupper_k}{1+\mupper_k} \; = \frac{2 \mupper_k}{1 + \mupper_k}
\end{align*}
Consequently, we have $\mupper_{k+1} = \frac{\mupper_k}{ (2
  \mupper_k)(1 + \mupper_k)} \; = \; \frac{1}{2} \big( 1 + \mupper_k
\big) = 1 + \frac{1}{2} \big(\mupper_k - 1 \big)$, as claimed.

From the decay rate~\eqref{EqnExpFast}, we see that taking $K = \lceil
\log \mupper_0 \rceil \leq 1 + \log_2 \mupper_0$ steps suffices to
ensure that $\mupper_K \leq 2$.  Since we also have $\Hess{k}(y)
\succeq \IdMat$ uniformly in $y$, it follows that $\sup_y
\mycond(\Hess{k}(y)) \leq \mupper_k$.  Consequently, this choice of
$K$ ensures that $\sup_y \mycond(\Hess{K}(y)) \leq 2$, as claimed.

\paragraph{Proof of the claim~\eqref{EqnBackwardControl}:}

By the representation~\eqref{EqnBackwardHess}
from~\Cref{LemForwardBackward} applied with $V = Y_{k+1}$ and $U =
Y_k$, we have $\BackHess_k(y)  = \Hess{k}(y) + \frac{a_k^2}{1 -
  a_k^2}$.  With our choice $a_k^2 = \frac{\mupper_k}{1 + \mupper_k}$,
we have $1 - a_k^2 = \frac{1}{1 + \mupper_k}$, and hence
$\frac{a_k^2}{1 - a_k^2} \; = \; \mupper_k$.  Since the eigenvalues of
$\Hess{k}(y)$ all lie in the interval $[1, \mupper_k]$ by
construction, the eigenvalues of $\BackHess_k(y)$ all lie in the
interval $[1 + \mupper_k, 2 \mupper_k]$, so that we are guaranteed to
have
\begin{align*}
  \mycond(\BackHess_k(y)) & \leq \frac{2 \mupper_k}{1 + \mupper_k} \leq
  2,
\end{align*}
as claimed.
\end{proof}


\subsection{Proof of~\Cref{LemSpectralControl}}
\label{SecProofLemSpectralControl}

\noindent We prove each of these two claims in turn.

\paragraph{Proof of the lower bound~\eqref{EqnVlow}:}

Combining the representation of $\BackHessU$ in
equation~\eqref{EqnBackwardHess} with our assumed lower bound on
$\HessU$, we have the uniform lower bound $\BackHessU(u) \succeq \big
(\mlower + \tfrac{a^2}{b^2} \big) \Id$.  Thus, the conditional
distribution $\pdensback$ is strongly log-concave, so that the
Brascamp--Lieb inequality (cf.~\Cref{SecBL}) is in force.  It allows
us to argue that
\begin{align*}
  \cov( \Uvar \mid \Vvar = v) & \stackrel{(i)}{\preceq}
  \Exs_{\pdensback}[ \big(\BackHessU(U) \big)^{-1} \big] \;
  \stackrel{(ii)}{\preceq} \; \big(\mlower + \tfrac{a^2}{b^2} \big)^{-1}
  \Id,
\end{align*}
where step (i) follows from the bound~\eqref{EqnSimpleBL}
in~\Cref{SecBL}; and step (ii) follows from the uniform lower bound on
$\BackHessU(u)$.

Consequently, we have the lower bound
\begin{align*}
  \Id - \frac{a^2}{b^2} \cov(\Uvar \mid \Vvar = v) & \succeq \Big(1 -
  \frac{a^2}{b^2} \frac{1}{\mlower + \frac{a^2}{b^2}} \Big) \IdMat \; =
  \; \Big( \frac{\mlower b^2 }{\mlower b^2 + a^2} \Big) \IdMat.
\end{align*}
Combining with the second-order Tweedie formula~\eqref{EqnForwardHess}
from~\Cref{LemForwardBackward}, we have shown that $\HessV(v) \succeq
\frac{\mlower}{ \mlower b^2 + a^2} \Id$, as claimed.

\paragraph{Proof of the upper bound~\eqref{EqnVup}:}

Returning to the second-order Tweedie
representation~\eqref{EqnForwardHess}, we see that upper bounds on
$\HessV$ require lower bounds on the covariance matrix $\cov(\Uvar
\mid \Vvar = v)$.  By a suitable embedding of our model into a
parametric family, we can obtain such a lower bound via the
Cramer--Rao approach, as we now describe.

Fix a realization $v \in \real^d$, and for each vector $\theta \in
\real^d$, define the shifted density $q_\theta(u) \defn \pdensback(u -
\theta)$.  (To keep the notation clean, we are suppressing the
dependence on $v$, since it remains fixed throughout the argument.)
We can now apply the Cramer--Rao bound to the parametric family $\{
q_\theta \mid \theta \in \real^d \}$.  By construction, we have
\begin{align*}
  -\nabla^2_\theta \log q_\theta(u) \Big|_{\theta = 0} & = - \nabla^2_u
  \log \pdensback(u) \; = \; \BackHessU(u).
\end{align*}
Consequently, the Fisher information for estimating $\theta = 0$ is
given by $\Fish \defn \Exs_{U \sim \pdensback} [\BackHessU(U)]$.  From
the representation~\eqref{EqnBackwardHess} of $\BackHessU$ and our
assumed upper bound on $\HessU$, we have
\begin{align*}
  \Fish & = \Exs_{U \sim \pdensback} \Big[ \HessU(U) + \frac{a^2}{b^2}
    \Id \Big] \; \preceq \Big( \mupper + \tfrac{a^2}{b^2} \Big) \Id,
\end{align*}
Inverting and negating the relation, we have $ -\Fish^{-1} \preceq -
\Big( \mupper + \tfrac{a^2}{b^2} \Big)^{-1} \Id$.

Now observe that $\psi(u) = u - \Exs[U \mid V= v]$ is an unbiased
estimate of $\theta = 0$ in this model, so that the Cramer--Rao bound
implies that $\cov(\Uvar \mid \Vvar = v) \succeq \Fish^{-1}$.  Putting
together the pieces, we have
\begin{align*}
  \Id - \tfrac{a^2}{b^2} \cov(\Uvar \mid \Vvar = v) \; \preceq \Id -
  \frac{a^2}{b^2} \Fish^{-1} \; \preceq \Big(1 - \frac{a^2}{b^2}
  \frac{1}{\mupper + \tfrac{a^2}{b^2}} \Big) \Id \; = \; \Big (
  \frac{\mupper b^2}{\mupper b^2 + a^2} \Big) \Id.
\end{align*}
Combining with the Tweedie form~\eqref{EqnForwardHess} of $\HessV$, we
conclude that $\HessV(v) \preceq \frac{\mupper}{\mupper b^2 + a^2}
\Id$, as claimed. \\



\subsection{Proof of~\Cref{LemConcaveError}}
\label{SecProofLemConcaveError}
We divide our proof into two parts, corresponding to each of the
two distances.

\paragraph{Proof for KL divergence:}

For the KL-divergence, let $Q$ be the joint distribution over $(Y_k,
Y_{k+1})$ defined by the marginal $q_{k+1}$ and the backward kernel
$\Qback_k$, so that $Y_k$ has marginal $q_k$.  Similarly, let $P$ be
the joint $(Y_k, Y_{k+1})$ defined by the marginal $p_{k+1}$ and the
backward kernel $\Pback_k$.  By the data-processing inequality, we
have $\kull(q_k \mid p_k) \leq \kull(Q \mid P)$.  Combined with the
chain rule for KL divergence, we find that
\begin{align*}
\kull(q_k \mid p_k) & \leq \Exs_{q_{k+1}} \Big[ \kull \big(
  \Qback_k(e_Y) \| \Pback_k(e_Y) \big) \Big] + \kull(q_{k+1}
\mid p_{k+1}) \\
& \leq \delta_k + \kull(q_{k+1} \mid p_{k+1}),
\end{align*}
as claimed.

\paragraph{Proof for Wasserstein $\Wass_2$:}
The Wasserstein distance satisfies the triangle inequality, so that we
can write
\begin{align*}  
  \Wass_2(q_k, p_k) = \Wass_2 \big( \Qback_k(q_{k+1}),
  \Pback_k(p_{k+1}) \big) & \leq \Wass_2 \big(\Qback_k(q_{k+1}),
  \Pback_k(q_{k+1}) \big) + \Wass_2 \big(\Pback_k(q_{k+1}),
  \Pback_k(p_{k+1}) \big).
\end{align*}
We have
\begin{align*}
\Wass_2^2( \Qback_k(q_{k+1}), \Pback_k(q_{k+1}) ) &
\stackrel{(i)}{\leq} \Exs_{q_{k+1}} \Wass_2^2 \big( \Qback_k(\cdot
\mid Y), \Pback_k( \cdot \mid Y) \big) \stackrel{(ii)}{\leq}
\delta_k^2
\end{align*}
where step (i) follows from joint convexity of $\Wass_2^2$ in its
arguments, and step (ii) follows from the assumed sampler
accuracy~\eqref{EqnSamplerAccuracy}.  Taking square roots of this
upper bound, and combining the two bounds yields
\begin{align}
\label{EqnWassIntermediate}  
  \Wass_2(q_k, p_k) & \leq \delta_k + \Wass_2
  \big(\Pback_k(q_{k+1}), \Pback_k(p_{k+1}) \big)
\end{align}

To complete the proof, it remains to show that $\Wass_2
\big(\Pback_k(q_{k+1}), \Pback_k(p_{k+1}) \big) \leq
\Wass_2(q_{k+1}, p_{k+1})$.  We make of~\Cref{LemWassStable} proved
in~\Cref{SecWassStable}, which controls the Wasserstein stability of
the kernels that arise in our backward analysis.  It allows us to
prove that
\begin{align}
  \label{EqnWassBound}
\Wass_2(\Pback_k(q_{k+1}), \Pback_k(p_{k+1})) & \leq a_k
\Wass_2(q_{k+1}, p_{k+1}) \; \leq \; \Wass_2(q_{k+1}, p_{k+1}),
\end{align}
and the remainder of the proof goes through as in the KL case.

In order to prove inequality~\eqref{EqnWassBound}, we first show how
our backward kernel can be converted to the form assumed
in~\Cref{LemWassStable}.  Let $\propto$ denote the proportionality
relation, keeping only terms dependent on $y_k$. Since $Y_{k+1} \mid
y_k \sim \Normal(a_k y_k, \IdMat/(1-a_k^2))$ by construction, we can
write
\begin{align*}
 p_{k \mid k+1}(y_k \mid y_{k+1}) & \propto \exp \Big \{ \log p_k(y_k)
 - \tfrac{1}{2 (1 - a_k^2)} \|a_k y_k - y_{k+1}\|_2^2 \Big \} \\
& \propto \exp \Big \{ -\psi(y_k) + \tfrac{a_k}{ (1 - a_k^2)}
 \inprod{y_k}{y_{k+1}} \Big \},
\end{align*}
where $\psi(y_k) = - \log p_k(y_k) + \frac{a_k^2}{2 (1 - a_k^2)}
\|y_k\|_2^2$.  Since $-\nabla^2 \log p_k(y_k) = \PlainHess_k(y_k)
\succeq \IdMat$, we have the lower bound $\nabla^2 \psi(y_k) \succeq 1
+ \frac{a_k^2}{1 - a_k^2}$, so that we can apply~\Cref{LemWassStable}
with $\alpha = 1 + \frac{a_k^2}{1- a_k^2}$ and $\beta = \frac{a_k}{1 -
  a_k^2} \geq 0$.  Finally, we observe that $\frac{\beta}{\alpha} =
\frac{ a_k/(1-a_k^2)}{1 + \tfrac{a_k^2}{1-a_k^2}} \; = \; a_k$, so
that the claim~\eqref{EqnWassBound} follows by application
of~\Cref{LemWassStable}.


\section{Auxiliary results related to~\Cref{ThmMultiModal}}

In this section, we collect the proofs of the auxiliary results
related to~\Cref{ThmMultiModal}.

\subsection{Proof of~\Cref{LemMultiForward}}
\label{SecProofLemMultiForward}

Define the sequence $\theta^2_0 = a_0^2 = 1/2$ and $\theta^2_k = a_k^2
\theta^2_{k-1}$ for $k = 1, 2, \ldots$.  We split our proof into two
parts.  By induction on $k$, it is easy to show that
\begin{align}
\label{EqnCanada}  
Y_{k} = \theta_{k-1} X + \sqrt{1 - \theta_{k-1}^2} W'_{k-1} \qquad \mbox{where $W'_{k-1}
  \sim \Normal(0, \IdMat)$ is independent of $X$.}
\end{align}
We make use of this representation repeatedly.

\paragraph{Proof of the forward sandwich~\eqref{EqnFhessian}:}

Applying the second-order Tweedie formula~\eqref{EqnForwardHess}
from~\Cref{LemForwardBackward} to the
representation~\eqref{EqnCanada}, we find that
\begin{align*}
\PlainHess_k(y) \;=\; \frac{1}{1- \theta_{k-1}^2}\Big\{ \IdMat \;-\;
\frac{\theta_{k-1}^2}{1 - \theta_{k-1}^2} \Cov\big(X \mid Y_k = y\big)
\Big\}.
\end{align*}
Since the covariance is positive semidefinite, it follows immediately
that $\PlainHess_k(y) \preceq \tfrac{1}{1 - \theta^2_{k-1}} \IdMat \;
\preceq \; 2 \IdMat$, where the final inequality follows since
$\theta^2_{k-1} \leq 1/2$ for all $k = 1, 2, \ldots$.

As for the lower bound, we have $1 \leq \frac{1}{(1 -
  \theta^2_{k-1})^2} \leq 4$, and hence
\begin{align*}
\PlainHess_k(y) \succeq \Big \{ 1 - 4 \theta^2_{k-1} \sup_{y \in
  \real^d} \opnorm{\cov(X \mid Y_k = y)} \Big \} \IdMat & = \; \big
\{1 - 4 \theta^2_{k-1} \CovCon_k \big \} \IdMat \; \succeq \;
\underbrace{-4 \theta^2_{k-1} \CovCon_k \IdMat}_{\equiv - \IterB_k},
\end{align*}
as claimed.

\paragraph{Proof of the backward sandwich~\eqref{EqnBhessian}:}

In this case, we use the equation $Y_{k+1} = a_k Y_k + \sqrt{1 -
  a_k^2} W_k$.  By the backward Hessian
representation~\eqref{EqnBackwardHess} from~\Cref{LemForwardBackward}
applied with $U = Y_k$, $V = Y_{k+1}$, $a^2 = a_k^2$ and $b^2 = 1 -
a_k^2$, we find that
\begin{align}
  \label{EqnBkExpand}
  \BackHess_k(y_k) \;=\; \PlainHess_k(y_k) \;+\; \frac{a_k^2}{1 -
    a_k^2} \IdMat \qquad \mbox{where $s_k \;\defn\; \frac{a_k^2}{1 -
      a_k^2}$.}
\end{align}
The sandwich condition~\eqref{EqnFhessian} ensures that $(-\IterB_k +
s_k)\, \IdMat \;\preceq\; \BackHess_k(y_k) \;\preceq\; (2 + s_k)\,
\IdMat$.  With the adaptively chosen stepsizes~\eqref{EqnMultiStepsize},
we have
\begin{align*}
  s_k \;=\; \tfrac{2\IterB_k + 2 }{(2 \IterB_k+2) + 1} \Big( 1 - \tfrac{2
    \IterB_k + 2 }{(2 \IterB_k+2) + 1} \Big)^{-1} \;=\; 2 \IterB_k + 2.
\end{align*}
Combining the pieces yields
\begin{align*}
\big( -\IterB_k + 2 \IterB_k +2 \big) \IdMat \;=\; \big(\IterB_k + 2) \IdMat
\;\preceq\; \BackHess_k(y_k) \;\preceq\; \big(2 + (2 \IterB_k + 2) \big)
\IdMat \; = \; 2 \, \big(\IterB_k + 2\big) \, \IdMat,
\end{align*}
as claimed.


\subsection{Proof of~Corollary\ref{CorMultiWorst}}
\label{SecProofCorMultiWorst}
Define the auxiliary sequence $\{\theta_k \}_{k=0}^\infty$ via
$\theta^2_0 = a_0^2 = 1/2$ and $\theta^2_k = a_k^2 \theta^2_{k-1}$ for
\mbox{$k = 1, 2, \ldots$.}  Introducing the shorthand $\Mtar = 8
\Cmax$, it suffices to specify a choice of $\Kfinal$ that ensures
$\theta^2_{K-1} \leq 1/\Mtar$, so that~\Cref{ThmMultiModal} can be
applied.  The remainder of our argument is devoted to showing that
\begin{align}
\label{EqnThetaTargetSimple}
\Kfinal & \leq 1 + 4 \Cmax + 3 \Big\lceil \log_2\!\Big(
\frac{\Mtar}{2} \Big) \Big \rceil \; \leq \; 7 \, \big(1 + \Cmax \big)
\end{align}
rounds are sufficient.

For our analysis, it is convenient to make use of the auxiliary
sequence defined by $u_0 = 2 \Cmax$, and $u_k \; \defn\; 4 \Cmax \,
\theta_k^2$.  Suppose that we can establish that
\begin{subequations}
\begin{align}
\label{EqnWashington}  
  u_{\Kfinal-1} \leq \utarget \defn 4\Cmax/\Mtar, \qquad \mbox{with
    $\Kfinal$ from equation~\eqref{EqnThetaTargetSimple}.}
\end{align}
It then follows that $\theta^2_{\Kfinal-1} = u_{\Kfinal-1}/(4 \Cmax)
\leq 1/\Mtar$, as desired. \\

Our proof of the bound~\eqref{EqnWashington} is based on the following
descent guarantee: for $k = 1, 2, \ldots$, we have
\begin{align}
\label{EqnDescent}
u_{k} - u_{k-1} & \leq - g(u_{k-1}) \quad \mbox{where} \quad g(s)
\defn \frac{s}{2 (s + 3/2)}.
\end{align}
\end{subequations}
We return to prove it momentarily; taking it as given for the moment,
let us prove the bound~\eqref{EqnWashington}.

\paragraph{Proof of the bound~\eqref{EqnWashington}:}

Introduce the shorthand $\Tfinal = \Kfinal - 1$.  Our goal is to
establish that $u_{\Tfinal} \leq \utarget = 4\Cmax / \Mtar$.  Since
$u_0 = 2 \Cmax$, we have the ratio $u_0/\utarget = \Mtar/2$, and we
analyze evolution of the iterates $\{u_k \}_{k \geq 0}$ as they move
through a sequence of $J \defn \lceil \log_2(\Mtar/2) \rceil$ epochs.
Each epoch is constructed so that the value $u_k$ drops by a factor of
$1/2$ as it transitions from the start to the end of the
interval. Concretely, for $j = 1, \ldots, J+1$, we define the
intervals
\begin{align*}
\MyInt{j} & \defn \Big( \underbrace{\frac{u_0}{2^j}}_{ \equiv
  \vup{j}}, \quad \underbrace{\frac{u_0}{2^{j-1}}}_{\equiv \vup{j-1}}
\Big],
\end{align*}
so that $\vup{0} = u_0$ and $\vup{J} = u_0/2^J \leq \utarget$.  At the
terminal round $\Kfinal$, we will ensure that $u_{\Tfinal} \in
\MyInt{J+1}$, so that $u_\Tfinal \leq \vup{J} \leq \utarget$.

The total number of steps $\Tfinal$ is defined by the sum $\Tfinal =
\sum_{j=1}^J \Nint{j}$ where $\Nint{j}$ is the number of steps $k$ for
which $u_k \in \MyInt{j}$.  In order to bound $\Tfinal$, we need to
bound $\Nint{j}$.  During epoch $j$, we need to reduce the value of
$u$ from $\vup{j-1}$ down to $\vup{j}$, so that the total decrease is
given by $\delup{j} \defn \vup{j-1} - \vup{j} = \vup{j}$.  Observe
that the function $g$ from the descent condition~\eqref{EqnDescent} is
strictly increasing.  Consequently, for values of $u \in \MyInt{j}$,
we have $g(u) \; \geq \; g(\vup{j}) = \frac{ \vup{j} }{ 2(\vup{j} +
  3/2 ) }$.  Using this fact, we have the bound
\begin{align*}
\Nint{j} \;\leq \; \frac{ \delup{j} }{ g(\vup{j}) } = \frac{\vup{j}}{
  \vup{j}/(2(\vup{j}+3/2)) } = 2\big( \vup{j} + 3/2 \big).
\end{align*}
Using this upper bound and summing over the epochs $j=1,\dots,J$
yields
\begin{align*}
\Tfinal \; = \; \sum_{j=1}^J \Nint{j} \; \leq \; \sum_{j=1}^J 2\big(\vup{j} +
3/2\big) = 2 \sum_{j=1}^J \vup{j} \;+\; 3 J.
\end{align*}
Since $\sum_{j=1}^J \vup{j} = u_0 (1 - 2^{-J})$ and $u_0 = 2\Cmax$, we
obtain
\begin{align*}
\Tfinal \;\le\; 4 \Cmax\big(1 - 2^{-J}\big) \; + \; 3 J \; \leq \; 4 \Cmax + 3 \lceil
\log_2(\Mtar/2) \rceil.
\end{align*}
Since $\Kfinal = \Tfinal + 1$, we have proved the
claim~\eqref{EqnThetaTargetSimple}.

\paragraph{Proof of the descent bound~\eqref{EqnDescent}:}
It remains to establish the descent condition.  From the definition of
$\IterB_k$ and the assumption that $\CovCon_k \leq \Cmax$, we have
\begin{align*}
\IterB_k & = 4 \theta^2_{k-1} \CovCon_k \; \leq 4 \theta^2_{k-1} \Cmax \; = \;
u_{k-1}.
\end{align*}
For $s > 0$, define the function $f(s) = \frac{2 s + 2}{2 s + 3}$, and
note that $a_k^2 = f(\IterB_k)$ by construction.  Since $f$ is an increasing function
and $\IterB_k \leq u_{k-1}$, we have
\begin{align*}
a_k^2 \; = f(\IterB_k) \; \leq f(u_{k-1}) \; = \;
\frac{2u_{k-1}+2}{2u_{k-1}+3} \;=\; 1 - \frac{1}{2(u_{k-1} + 3/2)}.
\end{align*}
Multiplying both sides by $u_{k-1}$ yields the one-step recursion
\begin{align*}
u_k & = a_k^2\, u_{k-1} \leq u_{k-1}\Big(1 -
\frac{1}{2(u_{k-1} + 3/2)}\Big),
\end{align*}
and re-arranging yields the claim~\eqref{EqnDescent}.


\section{Second-order Tweedie formula}
\label{AppProofLemForwardBackward}

In this section, we provide a proof of the following
lemma:
\mygraybox{
  \begin{lemma}[Forward and conditional Hessians]
    \label{LemForwardBackward}
    We have the second-order Tweedie formula
    \begin{subequations}    
      \begin{align}
        \label{EqnForwardHess}
\hspace*{-0.5in}        
        \mbox{\bf{Second-order Tweedie:}} \qquad
        \underbrace{\HessV(v)}_{-
    \nabla^2 \log \pdensv(v)} & = \frac{1}{b^2} \Big \{ \Id -
  \frac{a^2}{b^2} \cov( \Uvar \mid \Vvar = v) \Big \}.
      \end{align}
      Moreover, we have
      \begin{align}
        \label{EqnBackwardHess}
\hspace*{-1in} \mbox{\bf{Backward conditional Hessian:}} \qquad \quad
\underbrace{\BackHessU(u, v)}_{-\nabla^2_u \log \pdensback(u \mid v)}
= \HessU(u) + \frac{a^2}{b^2} \Id.
      \end{align}
    \end{subequations}
  \end{lemma}
}  
\noindent The second-order Tweedie formula~\eqref{EqnForwardHess} is a
known result (see the
papers~\cite{Efron11,chen2023improved,debortoli2022convergence,benton2024nearly}
for variants), whereas the backward conditional Hessian
formula~\eqref{EqnBackwardHess} follows directly from the structure of
the joint distribution.  For completeness, we provide a proof here.

We prove this lemma using a slightly more general result, which we
begin by stating.  Given a pair of random vectors $(\Uvar,\Vvar)$,
suppose that the integral representation
\begin{align}
\label{EqnMarginalV}
  \pdensv(v) & = \int_{\real^{d}} \pvcondu(v \mid u)\, \pdensu(u)\,du
\end{align}
of the marginal density admits sufficient regularity so that
differentiation under the integral is permitted.

\mygraybox{
\begin{lemma}[Hessian identity for marginal log densities]
  \label{LemLogMarginalHessian}
Under the above conditions, for every $v$ in the support of $\Vvar$,
we have the identity
\begin{align}
  \label{EqnLogMarginalHessian}
  \nabla_v^2 \log \pdensv(v) &= \Exs\big[ G(\Uvar,v) \mid \Vvar = v
    \big] \;+\; \cov\big( s(\Uvar,v) \mid \Vvar = v \big).
\end{align}
where $s(u,v) \defn \nabla_v \log p_{\Vvar \mid \Uvar}(v \mid u)$ is
the conditional score, and $G(u,v) \defn \nabla_v^2 \log p_{\Vvar \mid
  \Uvar}(v \mid u)$ is the conditional Hessian.
\end{lemma}
}
\noindent We first use this lemma to prove the two
claims~\eqref{EqnForwardHess} and~\eqref{EqnBackwardHess}
from~\Cref{LemForwardBackward}.  In~\Cref{SecTreadmill}, we return
to prove~\Cref{LemLogMarginalHessian}.

\subsection{Proof of equation~\eqref{EqnForwardHess}:}

We apply~\Cref{LemLogMarginalHessian} to our generative model $\Vvar =
a \Uvar + b \Wvar$, where $\Wvar$ is standard Gaussian.  By
definition, we have $\big( \Vvar \mid \Uvar = u \big) \sim \Normal(a
u, b^2 \IdMat)$, so that
\begin{align*}
s(u,v) & = \nabla_v \log \pvcondu(v \mid u) = \nabla_v \Big \{
-\frac{1}{2 b^2} \|v- a u\|_2^2 \Big \} \; = \; \frac{au - v}{b^2},
\quad \mbox{and} \\
G(u,v) & = \nabla^2 \log \pvcondu(v \mid u) \; = \; - \frac{1}{b^2} \IdMat.
\end{align*}
Thus, we have $\cov(s(U, v) \mid V = v) = \cov \Big( \frac{a U}{b^2}
\mid V = v) \; = \; \frac{a^2}{b^4} \cov(U \mid V= v)$.  Substituting
into~\eqref{EqnLogMarginalHessian} yields $\nabla^2 \log \pdensv(v) =
-\frac{1}{b^2} \IdMat + \frac{a^2}{b^4} \cov(U \mid V = v)$, and
re-arranging yields the claim~\eqref{EqnForwardHess}.


\subsection{Proof of equation~\eqref{EqnBackwardHess}:}

Since $\pucondv(u \mid v) = \pvcondu(v \mid u) \pdensu(u)/\pdensv(v)$,
we have
\begin{align*}
\BackHess(u,v) \equiv -\nabla^2_u \log \pucondv(u \mid v) & =
-\nabla^2_u \log \pdensu(u) - \nabla^2_u \log \pvcondu(v \mid u) \\
& = \HessU(u) + \nabla^2_u \Big \{ \frac{1}{2 b^2} \|a u - v\|_2^2
  \Big \} \; = \; \HessU(u) + \frac{a^2}{b^2} \IdMat,
\end{align*}
as claimed.

\subsection{Proof of~\Cref{LemLogMarginalHessian}}
\label{SecTreadmill}
By chain rule, we can compute $\nabla_v \log \pdensv(v) = \nabla
\pdensv(v)/\pdensv(v)$, and hence
\begin{subequations}
\begin{align}
  \nabla_v^2 \log \pdensv(v) & = \frac{1}{\pdensv(v)}\,\nabla_v^2
  \pdensv(v) \;-\; \frac{1}{\pdensv(v)^2}\, \big( \nabla_v \pdensv(v)
  \big) \big( \nabla_v \pdensv(v) \big)^\top \notag \\
\label{EqnTreadmill}
  & = \frac{1}{\pdensv(v)}\,\nabla_v^2 \pdensv(v) \;-\; \big( \nabla_v
  \log \pdensv(v) \big) \big( \nabla_v \log \pdensv(v) \big)^T.
\end{align}
Suppose that we can show that
\begin{align}
\label{EqnFisherOne}
  \nabla_v \log \pdensv(v) & = \Exs\big[ s(\Uvar,v) \mid \Vvar = v
    \big], \quad \mbox{and} \\
  \label{EqnFisherTwo}
  \nabla_v^2 \pdensv(v) & = \pdensv(v)\,
  \Exs\big[ G(\Uvar,v) + s(\Uvar,v)s(\Uvar,v)^\top \mid \Vvar = v
    \big].
\end{align}
\end{subequations}
Substituting these expressions into our
decomposition~\eqref{EqnTreadmill} then yields the
claim~\eqref{EqnLogMarginalHessian}.

\paragraph{Proof of the relation~\eqref{EqnFisherOne}:}
Beginning with the representation~\eqref{EqnMarginalV},
differentiating under the integral yields the relation $\nabla_v
\pdensv(v) \;=\; \int \nabla_v p_{\Vvar \mid \Uvar}(v \mid
u)\,p_\Uvar(u)\,du$.  Next we observe that
\begin{subequations}
\begin{align}
\label{EqnRice}
\nabla_v \pvcondu(v \mid u) & = p_{\Vvar \mid \Uvar}(v \mid
u)\,\nabla_v \log \pvcondu(v \mid u) = \pvcondu(v \mid u)\,s(u,v),
\qquad \mbox{and hence} \\
\label{EqnUsefulOne}
\nabla_v \pdensv(v) & = \int \pvcondu(v \mid u) s(u,v)\,
\pdensu(u)\,du.
\end{align}
\end{subequations}
Using Bayes’ rule, we can rewrite this as $\nabla_v \pdensv(v) =
\pdensv(v)\, \Exs\big[ s(\Uvar,v) \mid \Vvar = v \big]$, and dividing
by $\pdensv(v)$ yields the claim~\eqref{EqnFisherOne}.

\paragraph{Proof of the relation~\eqref{EqnFisherTwo}:}
Differentiating $\pdensv$ twice under the integral yields the
expression \mbox{$\nabla_v^2 \pdensv(v) \;=\; \int \nabla_v^2
  \pvcondu(v \mid u)\, \pdensu(u)\,du$.}  Moreover, differentiating
equation~\eqref{EqnUsefulOne} yields
\begin{align*}
  \nabla_v \pdensv(v) = \nabla_v\big( \pvcondu(v \mid u) s(u,v) \big)
  & = \nabla_v \pvcondu(v \mid u) s(u,v) + \pvcondu(v \mid u) \nabla_v
  s(u,v) \\
  & = \pvcondu(v \mid u) \Big \{ s(u,v) s(u,v)^T +  G(u,v) \Big \},
\end{align*}
where we have used the fact that $\nabla_v s(u,v) = G(u,v)$ and the
relation~\eqref{EqnRice}.

Putting together the pieces, we have
\begin{align*}
  \nabla_v^2 \pdensv(v) & = \int p_{\Vvar \mid \Uvar}(v \mid u)\,
  \Big( G(u,v) + s(u,v)s(u,v)^\top \Big)\, p_\Uvar(u)\,du \\ &=
  \pdensv(v)\, \Exs\big[ G(\Uvar,v) + s(\Uvar,v)s(\Uvar,v)^\top \mid
    \Vvar = v \big],
\end{align*}
again using Bayes’ rule.  This completes the proof of the
claim~\eqref{EqnFisherTwo}.


\section{Elementary consequence of Brascamp--Lieb inequality}
\label{SecBL}
Consider a density on $\real^d$ of the form $p(x) \propto
\exp(-\psi(x))$ where $\psi$ is twice differentiable and strictly
convex.  The Brascamp--Lieb inequality asserts that for any
differentiable function $f: \real^d \rightarrow \real$, we have
\begin{subequations}
  \begin{align}
\label{EqnFullBL}    
  \var_p(f(X)) & \leq \Exs_p \Big[ \inprod{\nabla f(X)}{ (\nabla^2
      \psi(X))^{-1} \nabla f(X)} \Big].      
  \end{align}
  Let us derive an elementary consequence that leads to an upper bound
  on the matrix $\cov_p(X)$.  For any given $v \in \real^d$, define
  the linear function $f_v(x) = \inprod{v}{x - \mu}$, where $\mu \defn
  \Exs[X]$.  Computing the derivative $\nabla f_v(x) = v$ and then
  applying inequality~\eqref{EqnFullBL} yields the upper bound
\begin{align*}
v^T \cov_p(X) v \; = \; \var_p (f_v(X)) \; \leq \; v^T \Exs_p \big[
  (\nabla^2 \psi(X))^{-1} \big] v.
\end{align*}
Since the choice of $v \in \real^d$ was arbitrary, it follows that
$\cov_p(X) \preceq \Exs_p\Big[ (\nabla^2 \psi(X))^{-1} \Big]$.  In
particular, when $\nabla^2 \psi(x) \succeq m \IdMat$ uniformly in $x$,
it follows that
\begin{align}
\label{EqnSimpleBL}
\cov_p(X) & \preceq \frac{1}{m} \IdMat.
\end{align}
\end{subequations}

\comment{
\section{Blah}

\mygraybox{
\begin{lemma}[Forward recursion]
  \label{LemForwardGauss}
  For the recursion~\eqref{EqnYrecursion}, we have the exact
  representation
  \begin{align}
    Y_{k+1} & = \theta_k X \, + \, \sqrt{1-\theta_k^2}\, W
  \end{align}
  where $W \sim \Normal(0, \IdMat)$ is Gaussian and independent of $X$.
\end{lemma}
}
\begin{proof}
We prove the claim by induction.  For $k = 0$, the claim follows from
our definition of $Y_0$.  Suppose that the claim holds at some round
$k \geq 1$, so that $Y_k = \theta_{k-1} X + \sqrt{ 1 - \theta_{k-1}^2}
W'$, for some Gaussian $W'$ independent of $X$.   We then have
\begin{align*}
  Y_{k+1} \; = \; a_k Y_k + \sqrt{1- a_k^2} W_k & = a_k \Big \{
  \theta_{k-1} X + \sqrt{1 - \theta^2_{k-1}} W' \Big \} + \sqrt{1-a_k^2} W \\
  & = \theta_k X + \underbrace{a_k \sqrt{1 - \theta^2_{k-1}} W' + \sqrt{1-a_k^2} W}_{\defn Z}
\end{align*}
where we have used the fact that $\theta_k = a_k \theta_{k-1}$.  It remains to show
that $\cov(Z) = \sqrt{1- \theta_k^2} \IdMat$.  Since $W$ and $W'$ are independent,
we  have
\begin{align*}
  \cov(Z) & = \Big(a_k^2 (1 - \theta^2_{k-1}) + 1 - a_k^2 \Big) \IdMat
  \; = \; \Big( 1 - \theta^2_{k} \Big) \IdMat,
\end{align*}
again using the fact that $a_K^2 \theta^2_{k-1} = \theta^2_k$.
\end{proof}
}

\section{Wasserstein stability of the backward kernel}
\label{SecWassStable}

Consider the Markov kernel $q \mapsto \Markov(q)(\cdot) \defn
\int_{\real^d} \pdensback(\cdot \mid v) q(v) dv$ where, for some
$\hackbeta \in \real$, the conditional density takes the form
\begin{subequations}
\begin{align}
\pdensback(u \mid v) & \propto \exp \Big \{ - \psi(u) + \hackbeta
\inprod{u}{v} \Big \},
\end{align}
and the Hessian $\nabla^2 \psi$ satisfies
\begin{align}
  \nabla^2 \psi(u) & \succeq \hackalpha \IdMat \qquad \mbox{uniformly
    over $u \in \real^d$.}
\end{align}
\end{subequations}

\mygraybox{
\begin{lemma}[Wasserstein stability of the backward kernel]
  \label{LemWassStable}
Under the above conditions, we have  
\begin{align}
    \label{Eq:KernelContractionFinal}
\Wass_2 \big( \Markov(q), \Markov(\qtil) \big) & \leq \frac{
  |\hackbeta|}{\hackalpha}\, \Wass_2(q, \qtil) \qquad \mbox{valid for
  any pair of distributions $q, \qtil$.}
\end{align}
\end{lemma}
}

\begin{proof}
The Wasserstein distance is defined via an infimum over couplings;
thus, by a standard ``gluing'' argument, it is sufficient to show that
\begin{subequations}
\begin{align}
\label{EqnGlue}    
\Wass_2 \big( \pdensback(\cdot \mid v), \; \pdensback(\cdot \mid
\vtil) \big) & \leq \frac{|\hackbeta|}{\hackalpha} \; \|v - \vtil\|_2
\qquad \mbox{for all $v, \vtil \in \real^d$.}
\end{align}
In order to do so, it is convenient to introduce the $d$-dimensional
exponential family
\begin{align*}
r_\theta(u) & \defn \exp \big \{ \inprod{\theta}{u} - \psi(u) -
A(\theta) \big \}
\end{align*}
where $A(\theta) = \log \int e^{\inprod{\theta}{u}- \psi(u)} du$ is
the log normalization constant.  By construction, for each $v \in
\real^d$, we have $\pdensback(u \mid v) \equiv r_{\theta(v)}(u)$ where
$\theta(v) \defn \hackbeta v$.  Consequently, in order to prove the
bound~\eqref{EqnGlue}, it suffices to show that
\begin{align}
\label{EqnGlueExp}  
\Wass_2 \big( r_\theta, r_{\thetatil} \big) & \leq
\frac{1}{\hackalpha} \|\theta - \thetatil\|_2 \qquad \mbox{for each
  $\theta, \thetatil$.}
\end{align}
\end{subequations}

Since $\nabla^2 \psi(u) \succeq \hackalpha \IdMat$ uniformly in $u \in
\real^d$, the density $r_{\thetatil}$ is
$\hackalpha$-strongly-log-concave.  Hence, Talagrand's $T_2$-bound
holds. (In particular, the Bakry--Emery criterion implies that
$r_{\thetatil}$ satisfies a log-Sobolev inequality (LSI) with
parameter $2/\hackalpha$, and the Otto--Villani translation from LSI
to $T_2$ then implies that $r_{\thetatil}$ satisfies Talagrand's
$T_2$-inequality with parameter $2/\hackalpha$.  See Corollary 7.3 and
Theorem 8.12 in the survey paper by~\cite{GozlanLeonard2010} for
details.) Consequently, we have
\begin{subequations}
\begin{align}
\label{EqnT2}    
\Wass^2_2(r_\theta, r_{\thetatil}) & \leq \frac{2}{\hackalpha}
\Kull(r_\theta \| r_{\thetatil}),
\end{align}
where $\Kull$ denotes the Kullback-Leibler divergence between
$r_\theta$ and $r_{\thetatil}$.  From the exponential family
structure~\citep{WaiJor08}, we have $\Kull(r_\theta \| r_{\thetatil})
= A(\thetatil) - A(\theta) - \inprod{\nabla A(\theta)}{\thetatil -
  \theta}$.  Taking one more derivative, we can write
\begin{align}
\label{EqnKLTwo}  
\Kull(r_\theta \| r_{\thetatil}) & = \frac{1}{2} (\thetatil -
\theta)^\top \nabla^2 A(\thetabar) (\thetatil - \theta),
\end{align}
where $\thetabar$ is some vector on the line joining $\theta$ and
$\thetatil$.  Thus, it remains to bound the Hessian $\nabla^2 A$.
Again using standard properties of exponential
families~\citep{WaiJor08}, we have $\nabla^2 A(\thetabar) =
\cov_{\thetabar} (U)$, where $\cov_{\thetabar}$ denotes the covariance
computed under the exponential family density $r_{\thetabar}$.  Since
$r_{\thetabar}$ is defined by a potential function that is
$\hackalpha$-strongly-convex, the Brascamp--Lieb inequality (see
equation~\eqref{EqnSimpleBL} in~\Cref{SecBL}) can be applied to assert
that
\begin{align}
\label{EqnNewBL}  
\nabla^2 A(\thetabar) = \cov_{\thetabar}(U) & \; \preceq \;
\frac{1}{\hackalpha} \IdMat \; \qquad \mbox{uniformly for all
  $\thetabar$.}
\end{align}
\end{subequations}

Applying the Brascamp--Lieb bound~\eqref{EqnNewBL} to
equation~\eqref{EqnKLTwo}, we find that the KL can be upper bounded as
$\KL(r_{\theta} \mid r_{\thetatil}) \leq \frac{1}{2 \hackalpha}
\|\theta - \thetatil\|_2^2$.  Combining with our
$T_2$-bound~\eqref{EqnT2}, we find that
\begin{align*}
\Wass^2_2(r_\theta, r_{\thetatil}) & \leq \frac{2}{\hackalpha}
\Kull(r_\theta \| r_{\thetatil}) \; \leq \; \frac{2}{\hackalpha} \;
\Biggr \{ \frac{1}{2 \hackalpha} \|\theta - \thetatil\|_2^2 \Biggr \}
\; = \; \frac{1}{\hackalpha^2} \|\theta - \thetatil\|_2^2.
\end{align*}
Taking square roots, we have proved the desired
bound~\eqref{EqnGlueExp}.
\end{proof}


\section{Results on score perturbations}

This section is devoted to results on score perturbations, beginning
in~\Cref{SecProofLemRobust} with the proof of the robustness guarantee
stated in~\Cref{LemRobust}. This proof makes use of an auxiliary
result of independent interest: \Cref{LemPerturb} on perturbed drift
functions in~\Cref{SecPerturb}.

\subsection{Proof of~\Cref{LemRobust}}
\label{SecProofLemRobust}

By definition of the approximate and exact updates, we have $\qtil_k =
\qtil_{k+1} \Qbacktil_k$ and $p_k = p_{k+1} \Pback_k$.  Thus, we can
write
\begin{align}
    \dtv(\qtil_k, p_k) = \dtv \left( \qtil _{k+1} \Qbacktil_k,\,
    p_{k+1} \Pback_k \right)
   & \stackrel{(i)}{\leq} \dtv (\qtil_{k+1} \Qbacktil_k, \qtil_{k+1}
    \Pback_k ) + \dtv(\qtil_{k+1} \Pback_k, \, p_{k+1} \Pback_k) \notag \\
\label{EqnPancake}    
  & \stackrel{(ii)}{\leq} \dtv (\qtil_{k+1} \Qbacktil_k, \qtil_{k+1}
    \Pback_k ) + \dtv(\qtil_{k+1} \, p_{k+1})
\end{align}
where step (i) follows from the triangle inequality; and step (ii)
follows since the TV distance is non-expansive under application of
the Markov kernel $\Pback_k$.  As for the first term in the
inequality, since both measures are mixtures over the same input law
$\qtil_{k+1}$, we have
\begin{subequations}
\begin{align}
\label{EqnSpain}  
\dtv(\qtil_{k+1} \Qbacktil_k,\, \qtil_{k+1}P_k) & \leq \mathbb
\Exs_{Y\sim \qtil_{k+1}} \dtv \big(\Qbacktil_k(\cdot\mid Y),
\Pback_k(\cdot\mid Y) )
\end{align}
By triangle inequality, we have
\begin{align}
\dtv\big(\Qbacktil_k(\cdot\mid y), \Pback_k(\cdot\mid y) \big) & \leq
\dtv (\Qbacktil_k(\cdot\mid y), \Qback_k(\cdot\mid y)) + \dtv \big(
\Qback_k(\cdot\mid y), \Pback_k(\cdot\mid y) \big) \notag \\
\label{EqnEcuador}
& \leq \varepsilon_k + \dtv \big(
\Qback_k(\cdot\mid y), \Pback_k(\cdot\mid y) \big),
\end{align}
where the second step uses our $\varepsilon_k$-accuracy assumption on
the sampler $\Qbacktil_k$.  Next, we have
\begin{align}
\dtv^2 \big( \Qback_k(\cdot\mid y), \Pback_k(\cdot\mid y) \big) &
\stackrel{(i)}{\leq} \frac{1}{2} \dkl \big( \Qback_k(\cdot\mid y) \|
\Pback_k(\cdot\mid y) \big) \notag \\
\label{EqnCosta}
& \stackrel{(ii)}{\leq} \frac{1}{4 \scon_k} I(\Qback_k(\cdot \mid y)
  \| \Pback_k(\cdot \mid y)),
\end{align}
\end{subequations}
where step (i) follows from Pinsker's inequality; and step (ii)
follows since any $\scon_k$-log-concave density satisfies a
log-Sobolev inequality with parameter $\scon_k$.  By definition, the
distribution $\Pback_y(\cdot \mid y)$ has score function $f(x) =
\score{k}(x) -\frac{a_k}{1-a_k^2}(a_kx - y)$, whereas $\Qback_y(\cdot
\mid y)$ is stationary for the SDE with drift function $\fhat(x) =
\scorehat{k}(x) -\frac{a_k}{1-a_k^2}(a_kx - y)$.  Thus, the drift
$\fhat$ is a perturbation of the score function $f$, so that we can
apply~\Cref{LemPerturb} to assert that $I(\Qback_k(\cdot \mid y) \|
\Pback_k(\cdot \mid y)) \leq \|\scorehat{k} - \score{k}
\|_{L^2(\Qback_k( \cdot \mid y))}^2$.

Combining with the bounds~\eqref{EqnEcuador} and~\eqref{EqnCosta},
we have shown that
\begin{align*}
  \dtv\big(\Qbacktil_k(\cdot\mid y), \Pback_k(\cdot\mid y) \big) &
  \leq \varepsilon_k + \frac{1}{2 \sqrt{\scon_k}} \|\scorehat{k} -
  \score{k} \|_{L^2(\Qback_k( \cdot \mid y))}.
\end{align*}
Taking expectations with respect to $\qtil_{k+1}$ and using the
bound~\eqref{EqnSpain}, we find that
\begin{align*}
  \dtv(\qtil_{k+1} \Qbacktil_k,\, \qtil_{k+1}P_k) & \leq \varepsilon_k +
  \frac{1}{2 \sqrt{\scon_k}} \Exs_{Y \sim \qtil_{k+1}} \|\scorehat{k}
  - \score{k} \|_{L^2(\Qback_k( \cdot \mid y))} \\
& \leq \varepsilon_k + \frac{1}{2 \sqrt{\scon_k}} \|\scorehat{k} -
  \score{k} \|_{L^2(\qbar_k)},
\end{align*}
where the final step follows by Jensen's inequality, and the
definition of $\qbar_k$.
Combining with inequality~\eqref{EqnPancake} completes the proof.


\subsection{Stationary drift perturbations}
\label{SecPerturb}
Let $\pi$ be a smooth density with score function $f$, and let
$\pihat$ be a smooth stationary distribution for the perturbed
diffusion $d X_t = \fhat(X_t) \,dt + \sqrt{2} \, dB_t$.

\mygraybox{
\begin{lemma}[Stationary drift perturbations]
\label{LemPerturb}
Under standard regularity conditions permitting integration by parts,
we have
\begin{align}
\label{EqnPerturb}
  I(\pihat \| \pi) & \defn \int \big\| \nabla \log
  \tfrac{\pihat(x)}{\pi(x)} \big\|_2^2 \; \pihat(x) dx \; \leq \;
  \|\fhat - f \|_{L^2(\pihat)}^2.
\end{align}
\end{lemma}
}
\begin{proof}
  Introduce the shorthand $e(x) \defn \fhat(x) - f(x)$ for the
  perturbation, and $r(x) \defn \log\frac{\pihat(x)}{\pi(x)}$ for the
  log density ratio.  Observe that $I(\pihat \| \pi) = \|\nabla
  r\|_{L^2(\pihat)}^2$ by definition.  Moreover, since $\pihat$ is
  stationary for the perturbed diffusion, it satisfies the stationary
  Fokker--Planck equation
\begin{align*}
  0 = - \nabla \cdot \{ \fhat \pihat \} + \Delta \pihat \; = \;
  -\nabla\cdot\{(f + e) \pihat\} +\Delta \pihat .
\end{align*}
Now we have
\begin{align*}
\nabla \pihat & = \pihat \nabla \log \pihat \; = \; \pihat \nabla
\big( \log \pi + r) \; = \;
\pihat f  + \pihat \nabla r,
\end{align*}
using the definition of $r$, and the fact that $f$ is the score
function of $\pi$.  Consequently, the FP equation can be rewritten as
$0 = \nabla \cdot(\pihat \nabla r) - \nabla\cdot(e \pihat)$.
Multiplying by $r$ and integrating by parts\footnote{ Indeed, for any
smooth vector field $A$, integration by parts gives $\int_{\real^d}
r\,\nabla\cdot A\,dx = -\int_{\real^d}\langle \nabla r,A\rangle\,dx$
provided the boundary term $\lim_{R\to\infty} \int_{\partial B_R}
r\,\langle A,n\rangle\,dS$ vanishes.  Applying this identity first
with $A=\pihat\nabla r$, and then with $A =e \pihat$ yields $0 = \int
r\,\nabla\cdot(\pihat\nabla r)\,dx - \int r\,\nabla\cdot(e\pihat)\,dx
= -\int \|\nabla r\|_2^2\,d\pihat + \int \langle e,\nabla
r\rangle\,d\pihat $.  } gives
\begin{align*}
  0 = -\int \|\nabla r\|_2^2\, d \pihat + \int \langle e, \nabla
  r\rangle\,d \pihat
\end{align*}
Re-arranging, we find that
\begin{align*}
  I(\pihat \|\pi) & = \int \|\nabla r\|_2^2\, d \pihat \; \leq |\int
  \langle e, \nabla r\rangle\,d \pihat| \; \leq \; I(\widehat
  \pi\|\pi)^{1/2} \left( \int \|e\|_2^2\,d\widehat \pi \right)^{1/2},
\end{align*}
where the last step follows by applying the Cauchy--Schwarz inequality
in $L^2(\pihat)$.
\end{proof}

\end{document}